\documentclass[sn-mathphys,Numbered]{sn-jnl}

\usepackage{graphicx}%
\usepackage{multirow}%
\usepackage{amsmath,amssymb,amsfonts}%
\usepackage{amsthm}%
\usepackage{mathrsfs}%
\usepackage[title]{appendix}%
\usepackage{xcolor}%
\usepackage{textcomp}%
\usepackage{manyfoot}%
\usepackage{booktabs}%
\usepackage{algorithm}%
\usepackage{algorithmicx}%
\usepackage{algpseudocode}%
\usepackage{listings}%

\usepackage{nicefrac}
\usepackage{array}
\usepackage{makecell}

\def\udots{\mathinner{\mskip1mu\raise1pt\vbox{\kern7pt\hbox{.}}\mskip2mu       
    \raise4pt\hbox{.}\mskip2mu\raise 7pt\hbox{.}\mskip1mu}}                     
\def\downbar#1{
\setbox10=\hbox{$#1$}
	    \dimen10=\ht10 \advance\dimen10 by 2.5pt
	    \ifdim \dimen10<15pt 
	       \advance\dimen10 by -0.5pt
	       \dimen11=\dimen10
	       \advance\dimen10 by 2.5pt 
	       \lower \dimen11
	    \else \lower \ht10 \fi
            \hbox {\hskip 1.5pt \vrule height \dimen10 depth \dp10}}
\def\upbar#1{
\setbox10=\hbox{$#1$}
	    \dimen10=\ht10 \advance\dimen10 by \dp10 \advance\dimen10 by 2.5pt
	    \ifdim \dimen10<15pt 
	       \advance\dimen10 by 2pt \fi
	    \raise 2.5pt \hbox {\hskip -1.5pt \vrule height \dimen10}}

\def\rz{\Bbb{R}}
\def\qz{\Bbb{Q}}

\def\nz{\Bbb{N}}

\def\cz{\Bbb{C}}
\def\gz{\Bbb{Z}}

\def\im{\text{\tt i}}

\DeclareMathOperator{\sign}{sign}
\DeclareMathOperator{\sgn}{sgn}
\DeclareMathOperator{\Ln}{Ln}
\DeclareMathOperator{\Arccos}{Arccos}
\DeclareMathOperator{\Arcsin}{Arcsin}
\newcolumntype{?}{!{\vrule width 1.5pt}}
\renewcommand{\thefootnote}{}

\newtheorem{theorem}{Theorem}
\newtheorem{corollary}{Corollary}

\newtheorem{lemma}{Lemma}

\begin{document}

\title[Multiscale matrix pencils for separable reconstruction problems]{Multiscale matrix pencils for separable reconstruction problems}


\author[1,2]{\fnm{Annie} \sur{Cuyt}}\email{annie.cuyt@uantwerpen.be}

\author*[2]{\fnm{Wen-shin} \sur{Lee}}\email{wen-shin.lee@stir.ac.uk}

\affil[1]{\orgdiv{Computational and Engineering Mathematics}, \orgname{University of Antwerp}, \orgaddress{\street{Middelheimlaan 1}, 
\postcode{2020}, \state{Antwerp}, \country{Belgium}}}

\affil*[2]{\orgdiv{Division of Computing Science and Mathematics}, \orgname{University of Stirling},\\ \orgaddress{\city{Stirling}, \postcode{FK9 4LA}, \state{Scotland}, \country{United Kingdom}}}



\abstract{
The nonlinear inverse problem of exponential data fitting is separable
since the fitting function is a linear combination of parameterized 
exponential functions, thus allowing to solve for the linear coefficients
separately from the nonlinear parameters. The matrix pencil method, which
reformulates the problem statement into a generalized eigenvalue problem
for the nonlinear parameters and a structured linear system for the linear
parameters, is generally considered as the more stable method to solve the
problem computationally.
In Section 2 the matrix pencil associated with the classical
complex exponential fitting or sparse interpolation problem is summarized
and the concepts of dilation and translation are introduced to obtain
matrix pencils at different scales.
\smallskip

Exponential analysis was earlier 
generalized to the use of several polynomial basis
functions and some operator eigenfunctions. However,
in most generalizations a computational scheme in terms of an eigenvalue problem
is lacking.
In the subsequent Sections 3--6 the matrix pencil formulation, including the 
dilation and translation paradigm, is generalized to more functions. 
Each of these periodic, polynomial or special function classes 
needs a tailored approach, where optimal use
is made of the properties of the parameterized elementary or special function 
used in the sparse interpolation problem under consideration. 
\smallskip

With each generalization a structured linear
matrix pencil is associated, immediately
leading to a computational scheme for the nonlinear and linear parameters,
respectively from a generalized eigenvalue
problem and one or more structured linear systems. 
\smallskip

Finally, in Section 7 we illustrate the new methods. 
}




\maketitle

\section{Introduction}

The nonlinear inverse problems of complex exponential analysis
\cite{Ka.Ma:81,Plonka:18}
and sparse 
polynomial interpolation \cite{Blahut:79,Be.Ti:det:88}
from 
uniformly sampled values 
can both be traced back to the exponential fitting method
of de Prony from the 18-th century \cite{dePr:ess:95,Hi:int:87}:
\begin{equation}f_j:= f(t_j)= \sum_{i=1}^n \alpha_i \exp(\phi_i t_j),  \alpha_i, \phi_i
\in \rz, \quad t_j=j\Delta \in \rz,  j=0, \ldots, 2n-1. \label{prony}\end{equation}
The French nobleman de Prony solved the
problem by obtaining the $n$ nonlinear parameters $\phi_i$ 
from the roots of a
polynomial and the $n$ coefficients $\alpha_i$ as the solution of a
Vandermonde structured linear system. Almost 200 years later this basic
fitting problem, that plays an important role 
\cite{Is.Vy:exp:99,Pe.Sc:exp:10} in many
computational science disciplines, engineering applications and digital
signal processing,
was reformulated in terms of a generalized eigenvalue
problem \cite{Hu.Sa:mat:90}. This reformulation, which is referred to as
the matrix pencil method, is generally the most reliable one when solving the
exponential analysis problem computationally.
\smallskip

It is the property
$$\exp(\phi_i t_{j+1}) = \exp(\phi_i \Delta) \exp(\phi_i t_j)$$
of the building blocks $\exp(\phi_i t)$
in \eqref{prony} that allows to split the
nonlinear interpolation problem \eqref{prony} into two numerical linear
algebra problems, namely the separate
computation of the nonlinear parameters $\phi_i$ from a generalized
eigenvalue problem on
the one hand and the linear coefficients $\alpha_i$ from a structured
linear system on the other.
\smallskip

Problem statement \eqref{prony} was partially generalized, 
to the use of non-standard polynomial bases such as
the Pochhammer basis and Chebyshev and Legendre polynomials
\cite{La.Sa:spa:95,Gi.La.ea:sym:04,Im.Ka.ea:spa:18,Po.Ta:spa:14,Pe.Pl.ea:rep:13}
and to the use of some eigenfunctions of linear operators
\cite{Pe.Pl:gen:13,Pl.St.ea:rec:19,St.Pl:gen:20}. Many of these
generalizations are unified in the algebraic framework described in
\cite{Ku.Ro.ea:lea:20}.
\smallskip

What is lacking in most of the generalizations above, is a reformulation in
terms of numerical linear algebra problems. In this paper we carry 
the generalized eigenvalue formulation of \eqref{prony}, 
so essentially the matrix pencil method, to
linear combinations of the trigonometric functions cosine, sine,
the hyperbolic cosine and sine functions,
the Chebyshev (1-st, 2-nd, 3-rd, 4-th kind) and spread polynomials,
the Gaussian function, the sinc and gamma function.
In addition, we introduce the paradigm of a selectable
dilation $\sigma$ and translation $\tau$ of the interpolation points, as used
in refinable function theory. All of the above functions namely satisfy a
property similar to
$$\exp(\phi_i t_{\tau + (j+1)\sigma}) = \exp(\phi_i t_\tau) \exp^\sigma(\phi_i
\Delta) \exp(\phi_i t_{j\sigma}),$$
which allows to separate the effect of the scale $\sigma$ and the shift
$\tau$ 
on the estimation of the parameters
$\phi_i$ and $\alpha_i$. This multiscale option will prove to be useful in
several situations, as further detailed in Section 2.2.
\smallskip

In each of the subsequent sections on the trigonometric and hyperbolic
functions, polynomial functions, the Gaussian distribution,
and some special functions,
a different approach is
required to express the nonlinear inverse problem 
\begin{equation}f_j = \sum_{i=1}^n \alpha_i g(\phi_i; t_j) , \qquad \alpha_i, \phi_i \in
\cz, \quad t_j \in \rz \label{genprony}\end{equation}
under consideration, as a generalized eigenvalue problem,
tailored to the particular properties of the building block $g(\phi_i; t)$
in use. 
The interpolant is always computed directly
from the evaluations $f_j$ 
where the $t_j$ follow 
some regular interpolation point pattern 
associated with the specific function $g(\phi_i; t)$.

\section{Exponential fitting}

We first lay out how the whole theory works for the exponential
problem, where $g(\phi_i; t) = \exp(\phi_i t)$.

\subsection{Scale and shift scheme}

By a combination of \cite{Hu.Sa:mat:90} and \cite{Cu.Le:how:20}
we obtain the following.
Let 
$f(t)$ be given by
\begin{equation}f(t) = \sum_{i=1}^n \alpha_i \exp(\phi_i t), \qquad \alpha_i, \phi_i \in\cz
\label{exp_mod}\end{equation}
and let us sample $f(t)$ at the equidistant points $t_j=j\Delta$ for
$j=0, 1, 2, \ldots$ with $\Delta \in \rz^+$, or more generally at
$t_{\tau+j\sigma}=(\tau+j\sigma)\Delta$ with
$\sigma \in \nz$ and $\tau \in \gz$, where
the frequency content in \eqref{exp_mod} is limited by
\cite{Ny:cer:28,Sh:com:49}
\begin{equation}|\Im(\phi_i)|\Delta < \pi, \qquad i=1, \ldots,n,
\label{nyquist}\end{equation}
with $\Im(\cdot)$ denoting the imaginary part.
More generally, $\sigma$ and $\tau$ can belong
to $\qz^+$ and $\qz$ respectively, as discussed in Section 2.5. 
The values $\sigma$ and $\tau$ are called 
the scaling factor and shift term respectively. 
We denote the collected samples by
$$f_{\tau+j\sigma}:= f(t_{\tau+j\sigma}), \qquad j=0, 1, 2, \ldots$$
From  $\exp(\phi_i t_{j+1}) = \exp(\phi_i\Delta) \exp(\phi_i t_j)$
we find that 
$$f_{j+1} = \sum_{i=1}^n \alpha_i \exp(\phi_i\Delta) \exp(\phi_i
j\Delta),$$
or more generally for
$\sigma \in\nz$ and $\tau \in \gz$ that
\begin{equation}f_{\tau+j\sigma} = \sum_{i=1}^n \alpha_i
\exp(\phi_i\tau\Delta)\exp(\phi_i j\sigma\Delta). \label{separate}\end{equation}
Hence we see that the scaling $\sigma$ and the shift $\tau$ are separated in a
natural way when evaluating
\eqref{exp_mod} at $t_{\tau+j\sigma}$, a property that
plays an important role in the sequel.
The freedom to choose $\sigma$ and $\tau$ 
when setting up the sampling scheme,
allows to stretch, shrink and translate the otherwise uniform progression
of sampling points dictated by the sampling step $\Delta$.
\smallskip

The aim is now to estimate the model order $n$, and the parameters
$\phi_1, \ldots, \phi_n$ and $\alpha_1, \ldots,
\alpha_n$ in \eqref{exp_mod} from samples $f_j$ at a selection of points
$t_j$. 

\subsection{Generalized eigenvalue formulation}

In this and the next subsection we assume for a moment that $n$ was 
determined.
With $n, \sigma \in \nz, \tau \in \gz$ we define
\begin{equation}{_\sigma^\tau}H_n := \begin{pmatrix}
f_\tau & f_{\tau+\sigma} & \cdots & f_{\tau+(n-1)\sigma} \\ 
f_{\tau+\sigma} & & & \\
\vdots & \udots & & \vdots \\
f_{\tau+(n-1)\sigma} & \cdots & & f_{\tau+(2n-2)\sigma}
\end{pmatrix}. \label{hankel}\end{equation}
It is well-known that the Hankel matrix ${_\sigma^\tau}H_n$ can be decomposed as
\begin{align}
\label{factor}
{_\sigma^\tau}H_n &= V_n \Lambda_n A_n V_n^T, \notag\\ 
V_n &= \begin{pmatrix}
1 & \cdots  & 1 \\
\exp(\phi_1\sigma\Delta) & \cdots & \exp(\phi_n\sigma\Delta) \\
\vdots & & \vdots \\
\exp(\phi_1(n-1)\sigma\Delta) & \cdots & \exp(\phi_n(n-1)\sigma\Delta)
\end{pmatrix}, \\ 
A_n &= \text{diag}(\alpha_1, \ldots, \alpha_n), \notag\\
\Lambda_n &= \text{diag}(\exp(\phi_1\tau\Delta), \ldots,
\exp(\phi_n\tau\Delta)). \notag
 \end{align}
This decomposition on the one hand translates \eqref{separate} and on the
other hand connects it to a generalized eigenvalue problem:
the values $\exp(\phi_i\sigma\Delta)$ 
can be retrieved \cite{Hu.Sa:mat:90} as the
generalized eigenvalues of the problem
\begin{equation}\left( {_\sigma^\sigma}H_n \right) v_i = \exp(\phi_i \sigma\Delta) 
\left( {_\sigma^0}H_n \right) v_i, \qquad i=1, \ldots,n, \label{gep}\end{equation}
where $v_i$ are the generalized right eigenvectors. 
Setting up this
generalized eigenvalue problem requires the $2n$ samples $f_{j\sigma},
j=0, \ldots, 2n-1$. A similar statement
holds for the values $\exp(\phi_i \tau\Delta)$ from the linear pencil
$({_\sigma^\tau}H_n, {_\sigma^0}H_n)$.
In \cite{dePr:ess:95,Hu.Sa:mat:90} the choices
$\sigma=1$ and $\tau=1$ are made 
and then, from the generalized eigenvalues $\exp(\phi_i\Delta)$,
the complex numbers $\phi_i$ can be retrieved uniquely
because of the restriction $|\Im(\phi_i)|\Delta<\pi$. 
\smallskip

Choosing $\sigma > 1$ offers a number of advantages though, among which:
\begin{itemize}
\item reconditioning \cite {Cu.Le:how:20,Pl.Wa.ea:det:18,Dihab}
of a possibly ill-conditioned problem statement,
\item superresolution \cite{Cu.Le:how:20,Cu.Ts.ea:fai:18}
in the case of clustered frequencies, 
\item validation \cite{Br.Cu.ea:vex:20}
of the exponential analysis output for $n$ and $\phi_i, i=1, \ldots, n$,
\item the possibility to parallellize the estimation of the parameters
$\phi_i$ \cite{Br.Cu.ea:vex:20}. 
\end{itemize}
\smallskip

With $\sigma>1$ the $\phi_i$ cannot necessarily be retrieved uniquely
from the generalized eigenvalues $\exp(\phi_i \sigma\Delta)$ since 
$\left| \Im (\phi_i) \right| \sigma\Delta$ may well be larger than $\pi$.
Let us indicate how to solve that problem which is called aliasing. 

\subsection{Vandermonde structured linear systems}

For chosen $\sigma$ and with
$\tau=0$, the $\alpha_i$ are computed from the interpolation conditions
\begin{equation}\sum_{i=1}^n \alpha_i \exp(\phi_i t_{j\sigma}) = f_{j\sigma}, 
\qquad j=0, \ldots, 2n-1, \qquad \sigma \in \nz, \label{vdm}\end{equation}
either by solving the system in the least squares sense, in the presence
of noise, or by solving a subset of $n$ interpolation
conditions in the case of noiseless samples. 
The samples of $f(t)$ occurring in
\eqref{vdm} are the same samples as the ones used to fill
the Hankel matrices in \eqref{gep} with.
Note that 
$$\exp(\phi_i t_{j\sigma})=\left( \exp(\phi_i\sigma\Delta) \right)^j,$$
and that for fixed $\sigma$
the coefficient matrix of \eqref{vdm} is therefore a transposed 
Vandermonde matrix with nodes $\exp(\phi_i \sigma \Delta)$.
In a noisy context the Hankel matrices in \eqref{gep}
can also be extended to rectangular $N \times \nu$
matrices with $N>\nu \ge n$ and the generalized eigenvalue
problem can be considered in a least squares sense \cite{Ch.Go:gen:06}.
In that case the index $j$ in \eqref{vdm} runs from 0 to $N+\nu-1$.
\smallskip

Next, for chosen nonzero $\tau$,
a shifted set of at least $n$ samples $f_{\tau+j\sigma}$ is interpreted as
\begin{equation}f_{\tau+j\sigma} = \sum_{i=1}^n \left( \alpha_i
\exp(\phi_i\tau\Delta) \right) \exp(\phi_i j\sigma\Delta), \qquad
j=k, \ldots, k+n-1, \qquad \tau \in \gz, \label{shift}\end{equation}
where $k \in \{0, 1, \ldots, n\}$ is fixed.
Note that \eqref{shift} is merely a shifted version of the original
problem \eqref{exp_mod}, where the effect
of the shift is pushed into the coefficients of \eqref{exp_mod}. The latter is
possible because of \eqref{separate}.
From \eqref{shift}, having the same (but maybe less) coefficient matrix
entries as \eqref{vdm}, we compute the unknown coefficients
$\alpha_i\exp(\phi_i\tau\Delta)$.
From $\alpha_i$ and $\alpha_i\exp(\phi_i\tau\Delta)$ we then obtain
$${\alpha_i\exp(\phi_i\tau\Delta) \over \alpha_i} = \exp(\phi_i\tau\Delta),$$
from which again the $\phi_i$ cannot necessarily be extracted unambiguously 
if $\tau>1$. But the following can be proved \cite{Cu.Le:how:20}. 
\smallskip

Denote $s_{i,\sigma}:= \sign \left( \Im \left( \Ln \left( 
\exp(\phi_i\sigma\Delta) \right) \right) \right)$ 
and $s_{i,\tau}:= \sign \left( \Im \left( \Ln \left( 
\exp(\phi_i\tau\Delta) \right) \right) \right)$, where 
$\Ln(\cdot)$
indicates the principal branch of the complex natural logarithm and 
$\left|\Im \left( 
\Ln \left( \exp(\phi_i\sigma\Delta) \right) \right) \right| \le \pi$.
If $\gcd(\sigma, \tau)=1$, then the sets
$$S_i=\left\{ {1 \over \sigma\Delta}\Ln \left( \exp(\phi_i\sigma\Delta) \right) 
+ {2 \pi \im \over \sigma\Delta} \ell, \
\ell=- s_{i,\sigma} \lfloor \sigma/2 \rfloor, \ldots, 
0, \ldots, s_{i,\sigma} (\lceil \sigma/2 \rceil -1) \right\},$$
$$T_i=\left\{ {1 \over \tau\Delta}\Ln \left( \exp(\phi_i\tau\Delta) \right)
+ {2 \pi \im \over \tau\Delta} \ell, \ 
\ell=- s_{i,\tau} \lfloor \tau/2 \rfloor, \ldots, 
0, \ldots, s_{i,\tau} (\lceil \tau/2 \rceil -1) \right\},$$
which contain all the possible arguments for $\phi_i$ in 
$\exp(\phi_i\sigma\Delta)$ from \eqref{gep} and in $\exp(\phi_i\tau\Delta)$
from \eqref{shift} respectively, 
have a unique intersection \cite{Cu.Le:how:20}. 
How to obtain this unique element in the intersection and identify
the $\phi_i$ is detailed in \cite{Cu.Le:how:20,Br.Cu.ea:vex:20}. 
Convenient choices for $\sigma$ and $\tau$ depend somewhat 
on the noise level and their selection is also discussed in 
\cite{Br.Cu.ea:vex:20}.
\smallskip

So at
this point the nonlinear parameters $\phi_i, i=1, \ldots, n$ and the linear 
$\alpha_i, i=1, \ldots, n$ in \eqref{exp_mod} are
computed through the solution of \eqref{gep} and \eqref{vdm}, and if $\sigma >
1$ also \eqref{shift}. Remains to discuss how to determine $n$.

\subsection{Determining the sparsity}

What can be said about the number of terms $n$ in \eqref{exp_mod}, which is
also called the sparsity? From \cite[p.~603]{He:app:74} and
\cite{Ka.Le:ear:03} we know for general $\sigma$ that
\begin{align}
\label{SVDn}
&\det {^0_\sigma}H_\nu = 0 \text{ only accidentally}, \qquad \nu<n, \notag\\
&\det {^0_\sigma}H_n \neq 0, \\
&\det {^0_\sigma}H_\nu = 0, \qquad  \nu > n. \notag
 \end{align}
The regularity of ${^0_\sigma}H_n$ persists for any value of $\Delta$ when
collecting the samples to fill the matrix with,
while an accidental singularity of ${^0_\sigma}H_\nu$ with $\nu<n$ only
occurs for an unfortunate choice of $\Delta$ that makes the determinant
zero.
A standard approach to make use of this statement is to compute a singular
value decomposition of the Hankel matrix ${^0_\sigma}H_\nu$ and this for
increasing values of $\nu$.
In the presence of noise and/or clustered eigenvalues, this technique
is not always reliable and we need to consider rather large values of
$\nu$ for a correct estimate of $n$ \cite{Cu.Ts.ea:fai:18} or turn our
attention to some validation add-on
\cite{Br.Cu.ea:vex:20}.
\smallskip

With $\sigma=1$ and in the absence of noise,
the exponential fitting problem can be solved from $2n$
samples for $\alpha_1, \ldots, \alpha_n$ and $\phi_1, \ldots, \phi_n$ and
at least one additional sample to confirm $n$. As pointed out already, it
may be worthwhile to take $\sigma>1$ and throw in at least
an additional $n$ values
$f_{\tau+j\sigma}$ to remedy the aliasing. Moreover, if $\max_{i=1, \ldots, n}
|\Im(\phi_i)|$ is quite
large, then $\Delta$ may become so small that collecting the 
samples $f_j$ becomes expensive and so it may be more feasible to work with a
larger sampling interval $\sigma\Delta$.

\subsection{Computational variants}

Besides having $\sigma \in \nz$ and $\tau \in \gz$, more general choices
are possible. An easy practical generalization is when
the scale factor and shift term are rational numbers $\sigma/\rho_1 \in
\qz^+$ and $\tau/\rho_2 \in \qz$ respectively, with 
$\sigma, \rho_1, \rho_2 \in \nz$ and $\tau
\in \gz$. In that case the condition $\gcd(\sigma, \tau)=1$ for $S_i$ and
$T_i$ to have a unique intersection, is replaced by
$\gcd(\overline\sigma, \overline\tau)=1$ where
$\sigma/\rho_1=\overline\sigma/\rho, \tau/\rho_2=\overline\tau/\rho$ with
$\rho=\text{lcm}(\rho_1, \rho_2)$.
\smallskip

We remark that, although the sparse interpolation problem can be solved
from the $2n$ samples $f_j, j=0, \ldots, 2n-1$ when $\sigma = 1$, we need
at least an additional $n$ samples at the shifted locations
$t_{\tau+j\sigma}, j=k, \ldots, k+n-1$ when $\sigma >1$. The former is
Prony's original problem statement in \cite{dePr:ess:95} and the latter is
the generalization presented in \cite{Cu.Le:how:20}. The factorisation
\eqref{factor} allows some alternative computational schemes, which may
deliver a better numerical accuracy but demand somewhat more samples.
\smallskip

First we remark that the
use of a shift $\tau$ can of course be replaced by the choice of a second 
scale factor $\tilde\sigma$ relatively prime with $\sigma$. But this
option requires the solution of two generalized eigenvalue problems of
which the generalized eigenvalues need to be matched in a combinatorial
step. Also, the
sampling scheme looks different and requires the $4n-1$ sampling points
$$\{t_{j\sigma}, 0 \le j \le 2n-1\} \cup \{t_{j\tilde\sigma}, 0 \le j \le
2n-1\}, \qquad \gcd(\sigma,\tilde\sigma)=1.$$
A better option is to set up the generalized eigenvalue problem
\begin{equation}{^\tau_\sigma}H_n v_i = \exp(\phi_i\tau\Delta) {^0_\sigma}H_n v_i, 
\qquad i=1, \ldots, n \label{variant}\end{equation}
which in a natural way connects each eigenvalue $\exp(\phi_i\tau\Delta)$,
bringing forth the set $T_i$, to its associated eigenvector $v_i$ 
bringing forth the set $S_i$. The latter is derived from the quotient
of any two consecutive entries in the vector ${^0_\sigma}H_n v_i$ which is
a scalar multiple of
$$\alpha_i(1, \exp(\phi_i\sigma\Delta), \ldots, 
\exp(\phi_i (n-1)\sigma\Delta))^T.$$
Such a scheme requires the $4n-2$ samples
$$\{t_{j\sigma}, 0 \le j \le 2n-2\} \cup \{t_{\tau+j\sigma}, 0 \le j \le
2n-2\}, \qquad \gcd(\sigma,\tau)=1.$$
Note that the generalized eigenvectors $v_i$ are actually insensitive to
the shift $\tau$: the eigenvectors of \eqref{gep} and \eqref{variant} are
identical. This is a remarkable fact that reappears in each of the
subsequent (sub)sections dealing with other choices for $g(\phi_i; t)$.
\smallskip

We now turn our attention to the identification of other families of
parameterized functions and patterns of sampling points. 
We distinguish between trigonometric and hyperbolic,
polynomial and other important functions.
Our focus here is on the derivation of the mathematical theory and not on
the practical aspects of the numerical computation.

\section{Trigonometric functions}

The generalized eigenvalue formulation \eqref{gep} incorporating the
scaling parameter $\sigma$, was generalized to
$g(\phi_i; t)=\cos(\phi_i t)$ in \cite{Gi.La.ea:sym:04} for integer
$\phi_i$ only. Here we
present a more elegant full generalization for $\cos(\phi_i t)$ 
including the use of a shift $\tau$ as in \eqref{shift} to restore
uniqueness of the solution if necessary. 
In addition we generalize the scale and shift approach to the functions 
sine, cosine hyperbolic and sine hyperbolic.

\subsection{Cosine function}

Let $g(\phi_i; t) = \cos(\phi_i t)$ with $\phi_i \in \rz$ where
\begin{equation}|\phi_i| \Delta < \pi, \qquad i=1, \ldots, n. \label{realnyq}\end{equation}
Since $\cos(\phi_i t) = \cos(-\phi_i t)$, we are only interested in 
the $|\phi_i|, i=1, \ldots, n$, disregarding the sign of each $\phi_i$.
With $t_j = j\Delta$ 
we still denote 
\begin{equation}f_{\tau+j\sigma}:= \sum_{i=1}^n \alpha_i \cos(\phi_i
(\tau+j\sigma)\Delta), \label{cprony}\end{equation}
and because of  
\begin{equation}{1 \over 2} \cos(\phi_i t_{j+1}) + {1 \over 2} \cos(\phi_i t_{j-1}) = 
\cos(\phi_i \Delta) \cos (\phi_i t_j)
\label{gerel_c}\end{equation}
we now also introduce for fixed chosen 
$\sigma$ and $\tau$,
\begin{equation}
\label{separate_c}
F_{\tau+j\sigma} := F(\sigma,\tau; t_j) = {1 \over 2} f_{\tau+j\sigma} + 
{1 \over 2} f_{\tau-j\sigma}, \\
= \sum_{i=1}^n \alpha_i \cos(\phi_i\tau\Delta)\cos(\phi_i j\sigma\Delta). 
\end{equation}
Relation \eqref{gerel_c} deals with the case $\sigma=1$ and $\tau=1$,
while the expression $F_{\tau+j\sigma}$ is a generalization of
\eqref{gerel_c} for general $\sigma$ and $\tau$.
Observe the achieved separation in \eqref{separate_c}
of the scaling $\sigma$ and the shift $\tau$.
We emphasize that $\sigma$ and
$\tau$ are fixed before defining the
$F(\sigma, \tau;t_j)$. Otherwise the index $j$ cannot be associated
uniquely with the value \nicefrac{1}{2} $(f_{\tau+j\sigma}+f_{\tau-j\sigma})$.
\smallskip

Besides the Hankel structured ${^\tau_\sigma}H_n$, 
we introduce the Toeplitz structured
$${^\tau_\sigma}T_n := \begin{pmatrix} f_\tau & f_{\tau-\sigma} &
\cdots & f_{\tau-(n-1)\sigma} \\
f_{\tau+\sigma} & & & \\
\vdots & & \ddots & \vdots \\ f_{\tau+(n-1)\sigma} & & \cdots & f_\tau
\end{pmatrix} ,$$
which is symmetric when $\tau=0$. 
Now consider the structured matrix
\begin{equation}{^\tau_\sigma}C_n := {1 \over 4} \left( {^\tau_\sigma}H_n \right) + 
{1 \over 4} \left( {^{\phantom{-}\tau}_{-\sigma}}H_n \right) + 
{1 \over 4} \left( {^\tau_\sigma}T_n \right) +
{1 \over 4} \left( {^{\phantom{-}\tau}_{-\sigma}}T_n \right), \label{Cmat}\end{equation}
where ${^{\phantom{-}\tau}_{-\sigma}}T_n = {^\tau_\sigma}T_n^T$.
When $\tau=0$, the first two matrices in the sum
coincide and the latter two do as well.
Note that working directly with the cosine function instead of expressing
it in terms of the exponential as $\cos x = (\exp(\im x) + \exp(-\im
x))/2$, reduces the size of the matrices involved in the pencil from $2n$ to $n$.
\smallskip 

\begin{theorem}
The matrix $^\tau_\sigma C_n$ factorizes as
\begin{align*}
{_\sigma^\tau}C_n &= W_n L_n A_n W_n^T, \\ 
W_n &= \begin{pmatrix}
1 & \cdots  & 1 \\
\cos(\phi_1\sigma\Delta) & \cdots & \cos(\phi_n\sigma\Delta) \\
\vdots & & \vdots \\
\cos(\phi_1(n-1)\sigma\Delta) & \cdots & \cos(\phi_n(n-1)\sigma\Delta) 
\end{pmatrix}, \\ 
A_n &= \text{\rm diag}(\alpha_1, \ldots, \alpha_n), \\
L_n &= \text{\rm diag}(\cos(\phi_1\tau\Delta), \ldots, \cos(\phi_n\tau\Delta)).
\end{align*}
\end{theorem}

\noindent{\it Proof. } The proof is a verification of the matrix product entry at position
$(k+1, \ell+1)$ for $k, \ell=0, \ldots, n-1$:
\begin{align*}
\text{\nicefrac{1}{4}} &f_{\tau+(k+\ell)\sigma} + \text{\nicefrac{1}{4}}
f_{\tau-(k+\ell)\sigma} + \text{ \nicefrac{1}{4}} f_{\tau+(k-\ell)\sigma} +
\text{\nicefrac{1}{4}} f_{\tau+(-k+\ell)\sigma} \\
&= \text{\nicefrac{1}{2}} \sum_{i=1}^n \alpha_i \cos(\phi_i \tau\Delta) 
\cos(\phi_i (k+\ell)\sigma\Delta)
+ \text{\nicefrac{1}{2}} \sum_{i=1}^n \alpha_i \cos(\phi_i \tau\Delta) 
\cos(\phi_i (k-\ell)\sigma\Delta) \\
&= \sum_{i=1}^n \alpha_i \cos(\phi_i\tau\Delta) \cos(\phi_i k\sigma\Delta)
\cos(\phi_i \ell\sigma\Delta). \hspace{3pt}\Box 
\end{align*}
This matrix factorization translates \eqref{separate_c} and opens the
door to the use of a generalized eigenvalue problem:
the cosine equivalent of \eqref{gep} becomes 
\begin{equation}\left( {^\sigma_\sigma}C_n \right) v_i = \cos(\phi_i \sigma\Delta) 
\left( {^0_\sigma}C_n \right) v_i, \qquad i=1, \ldots, n, \label{gep_c}\end{equation}
where $v_i$ are the generalized right eigenvectors. Setting up
\eqref{gep_c} takes $2n$ evaluations $f_{j\sigma}$, as in the exponential
case. Before turning our attention to the
extraction of the $\phi_i$ from the generalized eigenvalues
$\cos(\phi_i\sigma\Delta)$, we solve two structured linear systems of
interpolation conditions.
\smallskip

The coefficients $\alpha_i$ in \eqref{cprony} are computed from
$$\sum_{i=1}^n \alpha_i \cos(\phi_i j\sigma\Delta) = f_{j\sigma}, \qquad 
j=0, \ldots, 2n-1, \quad \sigma \in \nz.$$
Making use of \eqref{separate_c}, the coefficients
$\alpha_i\cos(\phi_i\tau\Delta)$ are obtained from the shifted
interpolation conditions
\begin{equation}\sum_{i=1}^n \left( \alpha_i \cos(\phi_i\tau\Delta) \right) 
\cos(\phi_i j\sigma\Delta) = F_{\tau+j\sigma}, \qquad 
j=k, \ldots, k+n-1, \quad \tau \in \gz, \label{shiftcos}\end{equation}
where $k \in \{0, 1, \ldots, n\}$ is fixed.
While for $\sigma=1$ the sparse interpolation problem
can be solved from
$2n$ samples taken at the points $t_j=j\Delta, j=0, \ldots, 2n-1$, 
for $\sigma > 1$ additional
samples are required at the shifted locations $t_{\tau \pm j\sigma}=
(\tau \pm j\sigma)\Delta$ in order to resolve the ambiguity that arises when 
extracting the nonlinear parameters $\phi_i$ from the
values $\cos(\phi_i \sigma\Delta)$. 
The quotient
$${\alpha_i \cos(\phi_i \tau\Delta) \over \alpha_i}, \qquad i=1, \ldots, n$$
delivers the values $\cos(\phi_i\tau\Delta), i=1, \ldots, n.$
Neither from $\cos(\phi_i\sigma\Delta)$ nor from $\cos(\phi_i\tau\Delta)$
the parameters $\phi_i$ can necessarily be extracted uniquely
when $\sigma>1$ and $\tau>1$. But the following result is proved in the
appendix. 
\smallskip

If $\gcd(\sigma,\tau)=1$, the sets
$$S_i=\left\{ 
{1 \over \sigma\Delta} \Arccos\left( \cos(\phi_i\sigma\Delta) \right) 
+ {2 \pi \over \sigma\Delta} \ell, \ 
\ell=-\lfloor \sigma/2 \rfloor, 
\ldots, 0, \ldots, \lceil \sigma/2 \rceil -1
\right\},$$
$$T_i=\left\{ 
{1 \over \tau\Delta} \Arccos\left( \cos(\phi_i\tau\Delta) \right) 
+ {2 \pi \over \tau\Delta} \ell, \ 
\ell=-\lfloor \tau/2 \rfloor, 
\ldots, 0, \ldots, \lceil \tau/2 \rceil -1 
\right\}$$
containing all the candidate arguments for $\phi_i$ in 
$\cos(\phi_i\sigma\Delta)$
and $\cos(\phi_i\tau\Delta)$ respectively,
have at most two elements in their intersection. Here $\Arccos(\cdot) 
\in [0,\pi]$
denotes the principal branch of the arccosine function. In case two
elements are found, then it suffices to extend \eqref{shiftcos} to
$$\sum_{i=1}^n \left( \alpha_i \cos(\phi_i(\sigma+\tau)\Delta) \right)
\cos(\phi_i j\sigma\Delta) = F_{(\sigma+\tau)+j\sigma}, \qquad
j=k, \ldots, k+n-1,$$
which only requires the additional sample
$f_{\tau+(k+n)\sigma}$ as $f_{\tau-(k+n-2)\sigma}$ is already available. 
From this extension,
$\cos(\phi_i(\sigma+\tau)\Delta)$ can be obtained in the same way as
$\cos(\phi_i\tau\Delta)$. As explained in the appendix,
only one of the two elements in the intersection
of $S_i$ and $T_i$ fits the computed $\cos(\phi_i(\sigma+\tau)\Delta)$
since $\gcd(\sigma, \tau)=1$ implies that also $\gcd(\sigma,\sigma+\tau) =
1 = \gcd(\tau,\sigma+\tau)$. 
\smallskip

So the unique identification of the $\phi_i$ 
can require $2n-1$ additional samples at the shifted locations $(\tau\pm
j\sigma)\Delta, j=0, \ldots n-1$ if the intersections $S_i \cap T_i$ are all 
singletons, or $2n$ additional samples, namely at 
$(\tau\pm j\sigma)\Delta, j=0, \ldots, n-1$ and $(\tau+n\sigma)\Delta$ 
if at least one of the intersections $S_i \cap T_i$ is not a singleton.
\smallskip

The factorization in Theorem 1 immediately allows to formulate the
following cosine analogue of \eqref{SVDn}.
\smallskip

\begin{corollary}
For the matrix ${^0_\sigma}C_n$ defined in \eqref{Cmat} holds that
$$\text{rank } {^0_\sigma}C_\nu = n, \qquad \nu \ge n.$$
\end{corollary}

To round up our discussion, we mention that from the factorization in
Theorem 1, it is clear that for the generalized eigenvector $v_i$ 
from the different generalized eigenvalue problem
$$\left({^\tau_\sigma}C_n\right) v_i = \cos(\phi_i\tau\Delta)
\left({^0_\sigma}C_n \right)v_i,$$
holds that ${^0_\sigma}C_n v_i$
is a scalar multiple of 
$$\alpha_i \left( 1, \cos(\phi_i\sigma\Delta),
\ldots, \cos(\phi_i(n-1)\sigma\Delta) \right)^T.$$
This immediately leads to a computational variant of the proposed scheme,
similar to the one given in Section 2.4 for the exponential function,
requiring somewhat more samples though.
Let us now turn our attention to other trigonometric functions.

\subsection{Sine function}

Let $g(\phi_i; t)=\sin(\phi_i t)$ and let \eqref{realnyq} hold. With
$t_j=j \Delta$ 
We denote
$$f_{\tau+j\sigma}:= \sum_{i=1}^n \alpha_i \sin(\phi_i
(\tau+j\sigma)\Delta),$$
and because of
\begin{equation}{1 \over 2} \sin(\phi_i t_{j+1}) + {1 \over 2} \sin(\phi_i t_{j-1}) =
\cos(\phi_i\Delta) \sin(\phi_i t_j), \qquad \Delta = t_{j+1}-t_j,
\label{gerel_s}\end{equation}
we introduce for fixed chosen $\sigma$ and $\tau$,
\begin{equation}
\label{separate_s}
F_{\tau+j\sigma} := F(\sigma, \tau; t_j)
={1 \over 2} f_{\tau+j\sigma} + {1 \over 2} f_{-\tau+j\sigma} \\
= \sum_{i=1}^n \left( \alpha_i \cos(\phi_i\tau\Delta) \right) \sin(\phi_i
j\sigma\Delta).
\end{equation} 
We fill the matrices ${_\sigma^\tau}H_n$ and the Toeplitz matrices
${_\sigma^\tau}T_n$ and define
\begin{equation}{_\sigma^\tau}B_n := {1 \over 4} \left(
{_{\phantom{+\tau}\sigma}^{\sigma+\tau}}H_n \right) + 
{1 \over 4} \left( {_{\phantom{-\tau}\sigma}^{\sigma-\tau}}H_n \right) + 
{1 \over 4} \left( {_{\phantom{+\tau}\sigma}^{\sigma+\tau}}T_n \right) + 
{1 \over 4} \left( {_{\phantom{-\tau}\sigma}^{\sigma-\tau}}T_n \right).
\label{Bmat}\end{equation}
\begin{theorem}
The structured matrix ${_\sigma^\tau}B_n$ factorizes as
\begin{align*}
{_\sigma^\tau}B_n &= U_n L_n A_n W_n^T, \\
U_n &= \begin{pmatrix} 
\sin(\phi_1\sigma\Delta) & \cdots & \sin(\phi_n\sigma\Delta) \\
\vdots & & \vdots \\
\sin(\phi_1 n\sigma\Delta) & \cdots & \sin(\phi_n n\sigma\Delta) 
\end{pmatrix} \\
W_n &= \begin{pmatrix}
1 & \cdots  & 1 \\
\cos(\phi_1\sigma\Delta) & \cdots & \cos(\phi_n\sigma\Delta) \\
\vdots & & \vdots \\
\cos(\phi_1(n-1)\sigma\Delta) & \cdots & \cos(\phi_n(n-1)\sigma\Delta)
\end{pmatrix}, \\ 
A_n &= \text{\rm diag}(\alpha_1, \ldots, \alpha_n), \\
L_n &= \text{\rm diag}(\cos(\phi_1\tau\Delta), \ldots,
\cos(\phi_n\tau\Delta)).
\end{align*} 
\end{theorem}

\noindent{\it Proof. } 
The proof is again a verification of the matrix product entry, at the
position $(k, \ell+1)$ with $k=1, \ldots, n$ and $\ell=0, \ldots, n-1$:
\begin{align*}
\text{\nicefrac{1}{4}} &f_{\tau+(k+\ell)\sigma} + \text{\nicefrac{1}{4}}
f_{-\tau+(k+\ell)\sigma} + \text{ \nicefrac{1}{4}} f_{\tau+(k-\ell)\sigma} +
\text{\nicefrac{1}{4}} f_{-\tau+(k-\ell)\sigma} \\
&= \text{\nicefrac{1}{2}} \sum_{i=1}^n \alpha_i \cos(\phi_i \tau\Delta) 
\sin(\phi_i (k+\ell)\sigma\Delta)
+ \text{\nicefrac{1}{2}} \sum_{i=1}^n \alpha_i \cos(\phi_i \tau\Delta) 
\sin(\phi_i (k-\ell)\sigma\Delta) \\
&= \sum_{i=1}^n \alpha_i \cos(\phi_i\tau\Delta) \sin(\phi_i k\sigma\Delta)
\cos(\phi_i \ell\sigma\Delta). \hspace{3pt}\Box 
\end{align*}

Note that the factorization involves precisely the building blocks 
in the shifted
evaluation \eqref{separate_s} of the help function $F(\sigma,\tau; t)$.
From this decomposition we find that the $\cos(\phi_i\sigma\Delta), i=1,
\ldots, n$ are obtained as the generalized eigenvalues of the problem
$$\left( {_\sigma^\sigma}B_n \right) v_i = \cos(\phi_i\sigma\Delta)
\left( {_\sigma^0}B_n \right) v_i, \qquad i=1, \ldots, n.$$ 
We point out that setting up this generalized eigenvalue problem requires
samples of $f(t)$ at the points $t_{(-n+1)\sigma}, \ldots, t_{2n\sigma}$. Since
$f(t_{j\sigma})=-f(t_{-j\sigma})$ and $f(0)=0$ it costs $2n$ samples.
Unfortunately, at this point we cannot compute the $\alpha_i, i=1, \ldots,
n$ from the linear system of interpolation conditions
$$\sum_{i=1}^n \alpha_i \sin(\phi_i j\sigma\Delta) = f_{j\sigma}, 
\qquad j=1, \ldots, 2n,$$
as we usually do,
because we do not have the matrix entries $\sin(\phi_i j\sigma\Delta)$ at
our disposal. It is however easy to obtain the values $\cos(\phi_i
j\sigma\Delta)$ 
because
$\cos(\phi_i j\sigma\Delta) = \cos\left( \pm j
\Arccos(\cos(\phi_i\sigma\Delta)) \right)$
where $\Arccos(\cos(\phi_i\sigma\Delta))$ returns the principal branch
value. The proper way to proceed is the following.
\smallskip

From Theorem 2 we get ${^0_\sigma}B_n^T = W_n A_n U_n^T$. So we
can obtain the $\alpha_i\sin(\phi_i\sigma\Delta)$ in the first column
of $A_n U_n^T$ from the structured linear system
\begin{equation}W_n \begin{pmatrix} \alpha_1\sin(\phi_1\sigma\Delta) \\ \vdots \\
\alpha_n\sin(\phi_n\sigma\Delta) \endpmatrix = \pmatrix b_{11} \\ \vdots
\\ b_{1n} \end{pmatrix}, \label{sincA}\end{equation}
where ${^0_\sigma}B_n= (b_{ij})_{i,j=1}^n$. From the generalized
eigenvalues
$\cos(\phi_i\sigma\Delta), i=1, \ldots$ and the
$\alpha_i\sin(\phi_i\sigma\Delta)$ we can now recursively compute for $j=1, \ldots, n,$
$$\alpha_i \sin(\phi_i j \sigma\Delta) = \alpha_i \sin(\phi_i
(j-1)\sigma\Delta) \cos(\phi_i\sigma\Delta) +
\cos(\phi_i(j-1)\sigma\Delta) \alpha_i\sin(\phi_i\sigma\Delta).$$ 
The system of shifted linear interpolation conditions
$$\sum_{i=1}^n \left( \alpha_i \cos(\phi_i\tau\Delta) \right)
\sin(\phi_i j\sigma\Delta) = F_{\tau+j\sigma}, 
\qquad j=k, \ldots, k+n-1, \quad 1 \le k \le n+1$$
can then be looked at as 
\begin{equation}\sum_{i=1}^n \left( \alpha_i \sin(\phi_i j\sigma\Delta) \right)
\cos(\phi_i\tau\Delta) = F_{\tau+j\sigma},
\qquad j=k, \ldots, k+n-1
\label{sincB}\end{equation}
having a coefficient matrix with entries $\alpha_i\sin(\phi_i j
\sigma\Delta)$ and unknowns $\cos(\phi_i\tau\Delta)$.
In order to retrieve
the $\phi_i$ uniquely from the values $\cos(\phi_i\sigma\Delta)$ and 
$\cos(\phi_i\tau\Delta)$ with $\gcd(\sigma,\tau)=1$, one proceeds as in the
cosine case. Finally, the $\alpha_i$ are obtained from the expressions
$\alpha_i\sin(\phi_i\sigma\Delta)$ after plugging in the correct arguments
$\phi_i$ in $\sin(\phi_i\sigma\Delta)$ and dividing by it. 
So compared to the previous sections, the intermediate 
computation of the $\alpha_i$ before knowing the $\phi_i$, is replaced by
the intermediate computation of the $\alpha_i\sin(\phi_i
\sigma\Delta)$. In the end, the $\alpha_i$ are revealed in a division,
without the need to solve an additional linear system.
\smallskip

From the factorization in Theorem 2, the following sine analogue of
\eqref{SVDn} follows immediately.
\smallskip

\begin{corollary}
For the matrix ${^0_\sigma}B_n$ defined in \eqref{Bmat} holds that
$$\text{rank } {^0_\sigma}B_\nu = n, \qquad \nu \ge n.$$
\end{corollary}

For completeness we mention that one also finds from this factorization 
that for $v_i$ in the generalized eigenvalue problem
$$\left( {^\tau_\sigma}B_n \right) v_i = \cos(\phi_i\tau\Delta) \left(
{^0_\sigma}B_n \right) v_i,$$
holds that ${^0_\sigma}B_n v_i$
is a scalar multiple of
$$\alpha_i \left( \sin(\phi_i\sigma\Delta), \ldots,
\sin(\phi_i n \sigma\Delta) \right)^T.$$

\subsection{Phase shifts in cosine and sine.}

It is possible to include phase shift parameters in the cosine and sine
interpolation schemes. We explain how, by working out the sparse
interpolation of 
\begin{equation}f(t) = \sum_{i=1}^n \alpha_i \sin(\phi_i t - \psi_i), \qquad 
\psi_i \in \rz, \quad \alpha_i, \phi_i \in \cz. \label{sinphase}\end{equation}
Since 
$\sin t = \left( \exp(\im t) - \exp(-\im t) \right)/2\im,$
we can write each term in \eqref{sinphase} as
$$\alpha_i \sin(\phi_i t-\psi_i) = {\alpha_i \exp(-\im\psi_i) \over 2\im}
\exp(\im \phi_i t) - {\alpha_i \exp(\im\psi_i) \over 2\im} \exp(-\im\phi_i
t).$$
So the sparse interpolation of \eqref{sinphase} can be solved by
considering the exponential sparse interpolation problem
$$\sum_{i=1}^{2n} \beta_i \exp(\im\zeta_i t),$$ 
where $\beta_{2i-1} = \alpha_i\exp(-\im\psi_i)/(2\im), \beta_{2i} =
-\alpha_i \exp(\im\psi_i)/(2\im)$ and $\zeta_{2i-1} = \phi_i=-\zeta_{2i}$.
The computation of the $\phi_i$ through the $\zeta_i$
remains separated from that of the
$\alpha_i$ and $\psi_i$. 
The latter are obtained as
\begin{align*}
\tan\psi_i &= -\im \frac{\beta_{2i} + \beta_{2i-1}}{\beta_{2i} -\beta_{2i-1}}, \\
\alpha_i &= -(\beta_{2i} + \beta_{2i-1})/\sin(\psi_i) = -\im (\beta_{2i} -
\beta_{2i-1})/\cos(\psi_i).
\end{align*}

\subsection{Hyperbolic functions}

For $g(\phi_i; t)=\cosh(\phi_i t)$ the computational scheme parallels that
of the cosine and for $g(\phi_i; t)=\sinh(\phi_i t)$ that of the sine. We
merely write down the main issues.
\smallskip

When $g(\phi_i; t)=\cosh(\phi_i t)$, let 
$$f_{\tau+j\sigma} := \sum_{i=1}^n \alpha_i \cosh(\phi_i
(\tau+j\sigma)\Delta)$$
and for fixed chosen $\sigma$ and $\tau$, let
$$
F_{\tau+j\sigma}:= {1 \over 2} f_{\tau+j\sigma} + {1 \over 2}
f_{\tau-j\sigma} 
= \sum_{i=1}^n \alpha_i \cosh(\phi_i\tau\Delta) \cosh(\phi_i
j\sigma\Delta).
$$
Subsequently the definition of the structured matrix ${_\sigma^\tau}C_n$ is
used and in the factorization of Theorem 1, the cosine function is
everywhere replaced by the cosine hyperbolic function.
\smallskip

Similarly, when $g(\phi_i; t)=\sinh(\phi_i t)$, let
$$f_{\tau+j\sigma} := \sum_{i=1}^n \alpha_i \sinh(\phi_i
(\tau+j\sigma)\Delta)$$
and for fixed chosen $\sigma$ and $\tau$, let
$$
F_{\tau+j\sigma}:= \frac{1}{2} f_{\tau+j\sigma} + \frac{1}{ 2}
f_{-\tau+j\sigma} 
= \sum_{i=1}^n \alpha_i \cosh(\phi_i\tau\Delta) \sinh(\phi_i
j\sigma\Delta).
$$
Now the definition of the structured matrix ${_\sigma^\tau}B_n$ is used
and in the factorization of Theorem 2 the occurrences of $\cos$ are
replaced by $\cosh$ and those of $\sin$ by $\sinh$.

\section{Polynomial functions}

The orthogonal Chebyshev polynomials were among the first polynomial
basis functions to be explored for use in combination with
a scaling factor $\sigma$, in the context of sparse interpolation in
symbolic-numeric computing
\cite{Gi.La.ea:sym:04}. We elaborate the topic further for numerical purposes
and for lacunary or supersparse interpolation, making use of the scale
factor $\sigma$ and the shift term $\tau$.
We also extend the approach to other
polynomial bases and connect to generalized eigenvalue formulations.

\subsection{Chebyshev 1st kind}

Let $g(m_i; t)=T_{m_i}(t)$ of degree $m_i$, which is defined by
$$T_m(t) = \cos(m\theta), \qquad t=\cos(\theta), \quad -1 \le t \le 1,$$
and consider the interpolation problem 
\begin{equation}f(t_j) = \sum_{i=1}^n \alpha_i T_{m_i}(t_j), \qquad \alpha_i \in \cz,
\quad m_i \in \nz. \label{T_mod}\end{equation}
The Chebyshev polynomials $T_m(t)$ satisfy the recurrence relation
$$T_{m+1}(t)= 2tT_m(t) - T_{m-1}(t), \qquad T_1(t) = t, \quad T_0(t) = 1$$
and the property
$${1 \over 2} T_{m_i}(t_{j+1}) + {1 \over 2} T_{m_i}(t_{j-1}) =
T_{m_i}(\cos\Delta) T_{m_i}(t_j).$$ 
With $0 \le m_1 < m_2 < \ldots < m_n <M$
we choose $t_j=\cos(j\Delta)$ where $0<\Delta\le\pi/M$. Note that the
points $t_j$ are allowed to occupy much more general positions than in 
\cite{Gi.La.ea:sym:04}. 
If $M$ is extremely
large and $n$ is small, in other words if the polynomial is very sparse, then
it is a good idea to recover the actual $m_i, i=1, \ldots, n$ in two
tiers as we explain now. Let 
$\gcd(\sigma,\tau)=1$. We denote
\begin{equation}f_{\tau+j\sigma}:= \sum_{i=1}^n \alpha_i T_{m_i}(t_{\tau+j\sigma})
\label{tcheb}\end{equation}
and introduce for fixed $\sigma$ and $\tau$,
$$
F_{\tau+j\sigma}:= F(\sigma, \tau; t_j) =
{1 \over 2} f_{\tau+j\sigma} + {1 \over 2} 
f_{\tau-j\sigma} 
= \sum_{i=1}^n \alpha_i T_{m_i}(\cos(\tau\Delta))
T_{m_i}(\cos(j\sigma\Delta)) $$
in order to separate the effect of $\sigma$ and $\tau$ in the evaluation. 
With the same matrices ${_\sigma^\tau}H_n, {_\sigma^\tau}T_n$ and
${_\sigma^\tau}C_n$ as in the cosine subsection, now filled with the
$f_{\tau+j\sigma}$ from \eqref{tcheb},
the values $T_{m_i}(\cos(\sigma\Delta))$ are the
generalized eigenvalues of the problem
\begin{equation}({_\sigma^\sigma}C_n) v_i = T_{m_i}(\cos(\sigma\Delta))
({_\sigma^0}C_n) v_i, \quad i=1, \ldots, n. \label{gep_T}\end{equation}
From the values $T_{m_i}(\cos(\sigma\Delta)) = \cos(m_i\sigma\Delta)$ the
integer $m_i$ cannot necessarily be retrieved unambiguously. 
We need to find
out which of the elements in the set
$$S_i= \left\{ {\pm 1 \over \sigma\Delta} \Arccos(\cos(m_i\sigma\Delta)) + 
{2\pi \over \sigma\Delta}\ell, \ 
\ell=0, \ldots, \sigma-1 \right\} \cap \gz_M$$
is the one satisfying \eqref{T_mod}, where
$\Arccos(\cos(m_i\sigma\Delta))/(\sigma\Delta) \le M/\sigma$.
Depending on the relationship between
$\sigma$ and $M$ (relatively prime, generator, divisor, $\ldots$) the
set $S_i$ may contain one or more candidate integers for $m_i$ evaluating
to the same value $\cos(m_i\sigma\Delta)$. 
To resolve the ambiguity we consider the 
Vandermonde-like system for the $\alpha_i, i=1, \ldots, n$,
$$\sum_{i=1}^n \alpha_i T_{m_i}(\cos(j\sigma\Delta)) =
f_{j\sigma}, \qquad j=0, \ldots, 2n-1,$$
and the shifted problem
$$\sum_{i=1}^n \left( \alpha_i T_{m_i}(\cos(\tau\Delta)) \right)
T_{m_i}(\cos(j\sigma\Delta)) = F_{\tau+j\sigma}, \qquad j=k, \ldots, k+n-1,
\quad \tau\in\gz,$$
from which we compute the $\alpha_i T_{m_i}(\cos(\tau\Delta)) = 
\alpha_i\cos(m_i\tau\Delta)$.
Then 
$$
\cos(m_i\tau\Delta) = T_{m_i}(\cos(\tau\Delta)) \\
= {\alpha_i T_{m_i}\cos(\tau\Delta) \over \alpha_i}, \qquad
i=1, \ldots, n.
$$
If the intersection of the set $S_i$ with the set
$$T_i= \left\{ {\pm 1 \over \tau\Delta} \Arccos(\cos(m_i\tau\Delta)) + 
{2\pi \over \tau\Delta}\ell, \ \ell=0, \ldots, \tau-1 \right\} \cap \gz_M$$
is processed as in Section 3.1, then one can 
eventually identify the correct $m_i$. An illustration thereof is given in
Section 7.3.
\smallskip

When replacing \eqref{gep_T} by
$$({_\sigma^\tau}C_n) v_i = T_{m_i}(\cos(\tau\Delta))
({_\sigma^0}C_n) v_i, \quad i=1, \ldots, n,$$
we find that for $v_i$ holds that ${^0_\sigma}C_n v_i$
is a scalar multiple of
$$\alpha_i \left( 1, T_{m_i}(\cos \sigma\Delta), \ldots, 
T_{m_i}(\cos (n-1)\sigma\Delta) \right)^T.$$
This offers an alternative algorithm similar to the alternative in
Section 3.1 on the cosine function.

\subsection{Chebyshev 2nd, 3rd and 4th kind}

While the Chebyshev polynomials $T_{m_i}(t)$ of the first kind are
intrinsically related to the cosine function, the Chebyshev polynomials
$U_{m_i}(t)$ of the second kind can be expressed using the sine function:
\begin{align*}
&U_m(t) = {\sin\left( (m+1)\theta \right) \over \sin\theta}, \qquad t =
\cos\theta, \quad -1 < t < 1,\\
&U_m(-1) = (-1)^m(m+1), \quad U_m(1)=m+1.
\end{align*}
Therefore the sparse interpolation problem
$$f(t_j) = \sum_{i=1}^n \alpha_i U_{m_i}(t_j), \qquad 
\alpha_i \in\cz, \quad m_i \in \nz$$
can be solved along the same lines as in Section 4.1 but now using the
samples 
$$f(t_{\tau+j\sigma}) \sin(\Arccos\; t_{\tau+j\sigma})$$ 
instead of the $f_{\tau+j\sigma}$, for the sparse interpolation of
$$\sum_{i=1}^n \alpha_i \sin \left( (m_i+1)\theta_j \right) = \sin\theta_j
f(t_j), \qquad t_j=\cos\theta_j.$$
In a very similar way, the sparse interpolation problems
$$f(t_j) = \sum_{i=1}^n \alpha_i V_{m_i} (t_j), \qquad \alpha_i \in \cz,
\quad m_i \in \nz,$$
$$f(t_j) = \sum_{i=1}^n \alpha_i W_{m_i} (t_j), \qquad \alpha_i \in \cz,
\quad m_i \in \nz$$
can be solved,
using the Chebyshev polynomials $V_{m_i}(t)$ and $W_{m_i}(t)$ of the 
third and fourth kind respectively, given by
\begin{align*}
V_m(t) &= \frac{\cos\left( (n+\text{\nicefrac{1}{2}}) \theta \right)}{\cos(\theta/2)}, \qquad t=\cos\theta, \quad -1< t \le 1, \\
W_m(t) &= \frac{\sin\left( (n+\text{\nicefrac{1}{2}}) \theta \right)}{\sin(\theta/2)}, \qquad t=\cos\theta, \quad -1\le t < 1. 
\end{align*}

\subsection{Spread polynomials}

Let $g(m_i;t)$ equal the
degree $m_i$ spread polynomial $S_{m_i}(t)$ on $[0,1]$, which is defined by
$$S_m(t) = \sin^2(m\theta), \qquad t=\sin^2(\theta), \quad 0 \le t \le
1.$$ 
The spread polynomials $S_m(t)$ are related to the Chebyshev polynomials of the
first kind by $1-2tS_m(t) = T_m(1-2t)$ and satisfy the recurrence relation
$$S_{m+1}(t) = 2(1-2t) S_m(t) - S_{m-1}(t) + 2t, \qquad S_1(t) = t, \quad
S_0(t) = 0$$
and the property
\begin{equation}S_m(t)S_r(t) = \text{\nicefrac{1}{2}} S_m(t) + \text{\nicefrac{1}{2}} S_r(t) -
\text{\nicefrac{1}{4}} S_{m+r}(t) - \text{\nicefrac{1}{4}} S_{m-r}(t).
\label{spreadgep}\end{equation}
We consider the interpolation problem
$$f(t_j) = \sum_{i=1}^n \alpha_i S_{m_i}(t_j), \qquad \alpha_i \in \cz,
\quad m_i \in \nz,$$
where $t_j=\sin^2(j\Delta), j=0, 1, 2, \ldots$ with $0 < \Delta \le \pi/(2M)$
and 
$0 < m_1 < \ldots < m_n < M.$
The 
$S_{m_i}(t) = \sin^2(m_i \arcsin\sqrt{t})$ satisfy
$${1 \over 2} (S_{m_i}(\sin^2 \Delta) + S_{m_i}(t_j)) - {1 \over 4}
(S_{m_i}(t_{j+1}) + S_{m_i}(t_{j-1})) = S_{m_i}(\sin^2\Delta) S_{m_i}(t_j).$$
As in Section 4.1 we present a two-tier approach, which for $\sigma\le 1$
reduces to one step and avoids the additional evaluations required for the
second step. However, as indicated above, the two-tier scheme offers some
additional possibilities.
We denote
$$f_{\tau+j\sigma} := \sum_{i=1}^n \alpha_i S_{m_i} \left(
\sin^2((\tau+j\sigma)\Delta) \right)
= \sum_{i=1}^n \alpha_i S_{m_i} \left(
S_{\tau+j\sigma}(\sin^2 \Delta) \right).$$
With 
$$F_{\tau+j\sigma} := F(\sigma, \tau; t_j)
= {1 \over 2} \left(f_\tau + f_{j\sigma}\right) -
{1 \over 4} \left( f_{\tau+j\sigma} + f_{\tau-j\sigma} \right)$$ 
we obtain
$$F_{\tau+j\sigma} = \sum_{i=1}^n \alpha_i S_{m_i}(\sin^2 \tau\Delta)
S_{m_i}(\sin^2 j\sigma\Delta).$$
So the effect 
of the scale factor $\sigma$ on the one hand and the shift term
$\tau$ on the other can again be separated in the evaluation
$F_{\tau+j\sigma}$.
\smallskip

We introduce the matrices 
\begin{equation}
\label{JKmat}
\begin{aligned}
{_\sigma}J_n &:= \left( \text{\nicefrac{1}{2}} f_{k\sigma} + \text{\nicefrac{1}{2}} 
f_{\ell\sigma} - \text{\nicefrac{1}{4}} f_{(k+\ell)\sigma} - \text{\nicefrac{1}{4}} 
f_{(k-\ell)\sigma} \right)_{k,\ell=1}^n, \\
{^\tau_\sigma}K_n &:= \left( \text{\nicefrac{1}{2}}F_{\tau+k\sigma} +
\text{\nicefrac{1}{2}}F_{\tau+\ell\sigma} - \text{\nicefrac{1}{4}}F_{\tau+(k+\ell)\sigma}
- \text{\nicefrac{1}{4}}F_{\tau+(k-\ell)\sigma}
\right)_{k,\ell=1}^n.
\end{aligned}
\end{equation}
\begin{theorem}
The matrices ${^\tau_\sigma}K_n$ and ${_\sigma}J_n$ factorize as
\begin{align*}
{^\tau_\sigma}K_n &= R_n L_n A_n R_n^T, \\
{_\sigma}J_n &= R_n A_n R_n^T, \\
R_n &= \begin{pmatrix} S_{m_1}(\sin^2\sigma\Delta) & \cdots &
S_{m_n}(\sin^2\sigma\Delta) \\ \vdots & & \vdots \\
S_{m_1}(\sin^2 n\sigma\Delta) & \cdots & 
S_{m_n}(\sin^2 n\sigma\Delta) \end{pmatrix}, \\ 
A_n &= \text{\rm diag} (\alpha_1, \ldots, \alpha_n), \\
L_n &= \text{\rm diag} \left( S_{m_1}(\sin^2 \tau\Delta), \ldots,
S_{m_n}(\sin^2 \tau\Delta) \right) \\
&= \text{\rm diag} \left( \sin^2(m_1\tau\Delta), \ldots,
\sin^2(m_n\tau\Delta) \right).
\end{align*}
\end{theorem}

\noindent{\it Proof. } 
The factorization is again verified at the level of the matrix entries, 
now  making use
of property \eqref{spreadgep}, which is slightly more particular. \hspace{3pt}$\Box$ 
\smallskip

This factorization paves the way to obtaining the values
$S_{m_i}(\sin^2 \sigma\Delta) = \sin^2(m_i\sigma\Delta)$ 
as the generalized eigenvalues of 
$$\left({^\sigma_\sigma}K_n\right) v_i = S_{m_i}(\sin^2 \sigma\Delta)
\left({_\sigma}J_n \right) v_i, \qquad i=1, \ldots, n.$$
Filling the matrices in this matrix pencil
requires $2n+1$ evaluations $f(j\sigma\Delta)$ for
$j=1 \ldots, 2n+1$. 
From these generalized eigenvalues we cannot necessarily uniquely 
deduce the values for the indices $m_i$.
Instead, we can obtain for each $i=1, \ldots, n$ the set of elements
\begin{multline*}
S_i = \left( \left\{ {\Arcsin(|\sin(m_i\sigma\Delta)|) \over \sigma\Delta} + 
{\pi \over \sigma\Delta}\ell, \ \ell=0, \ldots \lceil\sigma/2\rceil-1 \right\} 
\cup \right. \\ \left.
\left\{ {-\Arcsin(|\sin(m_i\sigma\Delta)|) \over \sigma\Delta} +
{\pi \over \sigma\Delta}\ell, \ \ell=1, \ldots \lfloor\sigma/2\rfloor \right\} 
\right) \cap \gz_M
\end{multline*}
characterising all the possible values for $m_i$
consistent with the sparse spread polynomial
interpolation problem. Fortunately, 
with $\gcd(\sigma, \tau)=1$, we can proceed as follows.
\smallskip

First, the coefficients $\alpha_i$ are obtained from the linear system of
interpolation conditions
$$\sum_{i=1}^n \alpha_i S_{m_i}(\sin^2 j\sigma\Delta) = f_{j\sigma},
\qquad j=0, \ldots, 2n-1.$$
The additional values $F_{\tau+j\sigma}$ lead to a second system of
interpolation conditions,
$$\sum_{i=1}^n \left( \alpha_i S_{m_i}(\sin^2 \tau\Delta) \right)
S_{m_i}(\sin^2 j\sigma\Delta) =  F_{\tau+j\sigma}, \qquad j=k, \ldots,
k+n-1, \quad 1 \le k,$$
which delivers the coefficients $\alpha_i S_{m_i}(\sin^2 \tau\Delta)$.
Dividing the two solution vectors of these linear systems componentwise
delivers the values $S_{m_i}(\sin^2 \tau\Delta), i=1, \ldots, n$ 
from which we obtain sets
\begin{multline*}
T_i = \left( \left\{ {\Arcsin(|\sin(m_i\tau\Delta)|) \over \tau\Delta} +
{\pi \over \tau\Delta}\ell, \ \ell=0, \ldots \lceil\tau/2\rceil-1 \right\}
\cup \right. \\ \left.
\left\{ {-\Arcsin(|\sin(m_i\tau\Delta)|) \over \tau\Delta} +
{\pi \over \tau\Delta}\ell, \ \ell=1, \ldots \lfloor\tau/2\rfloor \right\}
\right) \cap \gz_M
\end{multline*}
that have the correct $m_i$ in their intersection with the 
respective $S_i$. The proof of this statement follows a completely similar
course as that for the cosine building block $g(\phi_i;t)$, 
given in the Appendix. 
\smallskip

The factorization in Theorem 3 allows to write down 
a spread polynomial analogue of \eqref{SVDn}.
\smallskip

\begin{corollary}
For the matrices ${_\sigma}J_n$ and ${^0_\sigma}K_n$ defined in 
\eqref{JKmat} holds that
\begin{align*}
\text{rank } {_\sigma}J_\nu &= n, \qquad \nu \ge n, \\
\text{rank } {^0_\sigma}K_\nu &= n, \qquad \nu \ge n.
\end{align*}
\end{corollary}

To round up the discussion we mention that
from Theorem 3 and the generalized eigenvalue problem
$$\left({^\tau_\sigma}K_n\right) v_i = S_{m_i}(\sin^2 \tau\Delta)
\left({_\sigma}J_n \right) v_i, \qquad i=1, \ldots, n,$$
we also find that ${_\sigma}J_n v_i$
is a scalar multiple of
$$ \alpha_i \left( S_{m_i}(\sin^2\sigma\Delta),
\ldots, S_{m_i}(\sin^2 n\sigma\Delta) \right)^T.$$
At the expense of some additional samples this eigenvalue and
eigenvector combination offers
again an alternative computational scheme.

\section{Distribution functions}

In \cite[pp.~85--91]{Pe:gen:13} 
Prony's method is generalized from $g(\phi_i; t) = 
\exp(\phi_i t)$ 
with $\phi_i \in \cz$ to $g(\phi_i; t) = \exp(-(t-\phi_i)^2)$, 
to solve the interpolation problem
$$f(t_j) = \sum_{i=1}^n \alpha_i \exp\left( -{(t_j-\phi_i)^2 
\over 2w^2} \right), \qquad \alpha_i, \phi_i \in \cz,$$ 
with given fixed Gaussian peak width $w$.
Here we further generalize the algorithm to include the new scale and
shift paradigm. 
The scheme is useful
when modelling phenomena using Gaussian functions, as illustrated in
Section 7.1.
Without loss of generality we put $2w^2=1$. The easy adaptation to include
a fixed constant width factor in the formulas is left to the reader.
\smallskip

We again assume that \eqref{nyquist} holds, but now for $2\Delta$.
With $t_j= j\Delta$, the
Gaussian $g(\phi_i; t)=\exp(-(t-\phi_i)^2)$ satisfies
$$\exp \left( t_{j+1}^2 \right) \exp\left( -(t_{j+1}-\phi_i)^2 \right) =
\exp(2\phi_i\Delta)
\exp \left( t_j^2 \right) \exp\left( -(t_j-\phi_i)^2 \right).$$

Let us take a closer look at the evaluation of $f(t)$ at
$t_{\tau+j\sigma} = (\tau+j\sigma)\Delta, j=0, 1, \ldots$ with $\sigma \in
\nz$ and $\tau \in \gz$:
$$\exp\left( -((\tau+j\sigma)\Delta-\phi_i)^2 \right) = \exp \left(
-(\tau\Delta-\phi_i)^2 - j^2\sigma^2\Delta^2 -
2(\tau\Delta-\phi_i)j\sigma\Delta \right).$$
With the auxiliary function
\begin{equation}
\label{separate_d}
\begin{aligned}
F(\sigma, \tau; t_j) :&= \exp(2\tau j\sigma\Delta^2)
\exp(j^2\sigma^2\Delta^2) f(t_{\tau+j\sigma}) \\
&= \sum_{i=1}^n \left( \alpha_i \exp(-(\tau\Delta-\phi_i)^2) \right)
\exp(2\phi_i j\sigma\Delta)\, ,
\end{aligned}
\end{equation} 
we obtain a perfect separation of $\sigma$ and $\tau$ and 
the problem can be solved using
Prony's method. With fixed chosen $\sigma$ and $\tau$, the value
$F(\sigma,\tau; t_j)$ is denoted by $F_{\tau+j\sigma}$.
\smallskip

\begin{theorem}
The Hankel structured matrix
$${_\sigma^\tau}G_n := \begin{pmatrix}
F_\tau & F_{\tau+\sigma} & \cdots & F_{\tau+(n-1)\sigma} \\
F_{\tau+\sigma} & & & \\
\vdots & \udots & & \vdots \\
F_{\tau+(n-1)\sigma} & \cdots & & F_{\tau+(2n-2)\sigma}
\end{pmatrix}$$
factorizes as
\begin{align*} {_\sigma^\tau}G_n &= E_n L_n A_n E_n^T, \\
E_n &=  \begin{pmatrix} 1 & \cdots & 1 \\
\exp(2\phi_1 \sigma \Delta) & \cdots & \exp(2\phi_n \sigma \Delta) \\
\vdots &  & \vdots \\
\exp(2\phi_1 (n-1)\sigma \Delta) & \cdots & \exp(2\phi_n (n-1)\sigma \Delta) 
\end{pmatrix} \\ 
A_n &= \text{\rm diag} (\alpha_1\exp(-\phi_1^2), \ldots,
\alpha_n\exp(-\phi_n^2)) \\
L_n &= \text{\rm diag} \left(\exp(-\tau^2\Delta^2+2\tau\Delta\phi_1), \ldots,
\exp(-\tau^2\Delta^2+2\tau\Delta\phi_n) \right).
\end{align*}
\end{theorem}

\noindent{\it Proof. } 
The proof is again by verification of the entry $F_{\tau+(k+\ell)\sigma}$
in ${_\sigma^\tau G_n}$ at position $(k+1, \ell+1)$ for $k=0, \ldots, n-1$
and $\ell=0, \ldots, n-1$. \hspace{3pt}$\Box$ 
\smallskip

With $\tau=0,\sigma$ the values $\exp(2\phi_i\sigma\Delta)$ 
are retrieved as a factor of the generalized eigenvalues of the problem
$$\left( {_\sigma^\sigma}G_n \right) v_i = \exp(-\sigma^2\Delta^2)
\exp(2\phi_i\sigma\Delta) \left( {_\sigma^0}G_n\right) v_i, 
\qquad i=1, \ldots, n.$$
As we know from the exponential case, the $\phi_i$ cannot 
necessarily be identified
unambiguously from $\exp(2\phi_i\sigma\Delta)$ when $\sigma > 1$. 
In order to remedy that, 
we turn our attention to two structured linear systems. The
first one, where $\tau=0$, 
$$\sum_{i=1}^n \alpha_i \exp\left( -(t_{j\sigma}-\phi_i)^2 \right) =
f_{j\sigma}, \qquad j=0, \ldots, 2n-1,$$
delivers the $\alpha_i \exp(-\phi_i^2)$
after rewriting it as 
$$\sum_{i=1}^n \left( \alpha_i\exp(-\phi_i^2) \right) \exp(2\phi_i
j\sigma\Delta) = \exp(j^2\sigma^2\Delta^2) f_{j\sigma}= F(\sigma, 0; t_j),
\qquad j=0, \ldots, 2n-1.$$
The coefficient matrix of this linear system is Vandermonde structured
with entry $\left( \exp(2 \phi_i \sigma\Delta) \right)^j$ at position
$(j+1, i)$. 
The second linear system, where $\tau >0$,
delivers the $\alpha_i \exp(-(\tau\Delta-\phi_i)^2)$
through \eqref{separate_d},
$$\sum_{i=1}^n \left( \alpha_i \exp(-(\tau\Delta-\phi_i)^2) \right) 
\exp(2\phi_i j\sigma\Delta) = 
F_{\tau+j\sigma}, 
\qquad j=k, \ldots, k+n-1.$$
Here the coefficient matrix is structured identically as in the
first linear system.
From both solutions we obtain
\begin{align*}
\exp(\tau^2\Delta^2) &{\alpha_i \exp(-(\tau\Delta-\phi_i)^2) \over 
\alpha_i \exp(-\phi_i^2)} \\
&= \exp(\tau^2\Delta^2)
{\alpha_i\exp(-\tau^2\Delta^2)\exp(-\phi_i^2)\exp(2\phi_i\tau\Delta) \over
\alpha_i \exp(-\phi_i^2)} = \exp(2\phi_i\tau\Delta).
\end{align*}
From the values $\exp(2\phi_i\sigma\Delta), i =1, \ldots, n$ and
$\exp(2\phi_i\tau\Delta), i=1, \ldots, n$ the parameters $2\phi_i$ can be
extracted as explained in Section 2, under the condition that
$\gcd(\sigma,\tau)=1$. 
\smallskip

The values $\exp(2\phi_i\tau\Delta)$ and
$\exp(2\phi_i\sigma\Delta)$ can also be retrieved respectively from the
generalized eigenvalues and the generalized eigenvectors of the
alternative problem
$$\left( {^\tau_\sigma}G_n \right) v_i = \exp(-\tau^2\Delta^2)
\exp(2\phi_i\tau\Delta) \left( {_\sigma^0}G_n\right) v_i, 
\qquad i=1, \ldots, n,$$
with ${^0_\sigma}G_n v_i$ being a scalar multiple of
$$\alpha_i \left( 1, \exp(2\phi_i\sigma\Delta),
\ldots, \exp(2\phi_i(n-1)\sigma\Delta) \right)^T,$$
thereby requiring at least of $4n-2$ samples instead of $3n$ samples.
To conclude, the following analogue of \eqref{SVDn} can be given.
\smallskip

\begin{corollary}
For the matrix ${^0_\sigma}G_n$ given in Theorem 4 holds that
$$\text{rank } {^0_\sigma}G_\nu = n, \qquad \nu \ge n.$$
\end{corollary}

\section{Some special functions}

The sinc function is widely used in digital signal processing, especially
in seismic data processing where it is a natural interpolant.
There are 
several similarities between the narrowing sinc function and the Dirac delta 
function, among which the shape of the pulse. 
A large number of papers, among which \cite{Batenkov}, 
already discuss the determination of a so-called train of Dirac spikes 
and their amplitudes, which is essentially an exponential fitting problem.
This is generalized here to the use of the sinc function, including a
matrix pencil formulation.
\smallskip

The gamma function first arose in connection with the interpolation problem 
of finding a function that equals $n!$ when the argument is a positive
integer. Nowadays the
function plays an important role in mathematics, physics and engineering.
In \cite{La.Sa:spa:95} sparse interpolation or exponential analysis was
already generalized to the Pochhammer basis $(t)_m = t (t+1) \cdots
(t+m-1)$, also called rising factorial,
which is related to the gamma function by
$(t)_m = \Gamma(t+m)/\Gamma(t)$ for $t \not\in\gz^- \cup \{0\}$. 
Here we generalize
the method to the direct use of the gamma function and we present a matrix
pencil formulation as well.

\subsection{The sampling function $\sin(x)/x$}

Let $g(\phi_i; t)=\text{sinc}(\phi_i t)$ where $\text{sinc}(t)$ is
historically defined
by $\text{sinc}(t)=\sin t / t$. So our sparse interpolation problem is
$$f(t_j)=\sum_{i=1}^n \alpha_i \text{sinc}(\phi_i t_j), \qquad t_j =
j\Delta,$$
with the same assumptions for $\phi_i$ and $\Delta$ as in Section 3.
In order to solve this inverse problem of identifying the $\phi_i$ and
$\alpha_i$ for $i=1, \ldots, n$, we introduce
$$
F(t_j) := j\Delta f(t_j) 
= \sum_{i=1}^n \left( {\alpha_i \over \phi_i}\right) \sin(\phi_i j\Delta)
$$
and apply the technique from Section 3.2 for the 
separate identification of the
nonlinear parameters $\phi_i$ and linear parameters $\alpha_i/\phi_i$
in the sparse sine interpolation.

\subsection{The gamma function $\Gamma(z)$}

With the new tools obtained
so far, it is also possible to extend the theory to other functions such
as the gamma function $\Gamma(z)$.
The function $g(\phi_i; z) = \Gamma(z+\phi_i)$ with $z, \phi_i \in \cz$, 
satisfies the relation
\begin{equation}\Gamma(\Delta+1+\phi_i) = (\Delta +\phi_i) \Gamma(\Delta+\phi_i), \qquad 
\Delta \in \cz,  \Delta+\phi_i \in \cz\setminus \{0, -1, -2,
\ldots\}. \label{recgam}\end{equation}
Our interest is in the sparse interpolation of
$$f(z) = \sum_{i=1}^n \alpha_i \Gamma(z+\phi_i), \qquad z+\phi_i \in
\cz\setminus \{0, -1, -2, \ldots\}$$
where the $\alpha_i, \phi_i, i=1, \ldots, n$ are unknown.
In the sample point $z=\Delta$ we define
\begin{align*}
F_0(\Delta) :&= f(\Delta), \\
F_j(\Delta) :&= 
F_{j-1}(\Delta+1) - \Delta F_{j-1}(\Delta), 
\qquad j=1, 2, \ldots \end{align*}
If by the choice of $\Delta$, one or more of the $\Delta+\phi_i, i=1, \ldots, n$
accidentally belong to the set of nonpositive integers, then one cannot
sample $f(z)$ at $z=\Delta$. In that case a complex
shift $\tau$ can help out. It suffices to shift the arguments
$\Delta+\phi_i$ away from the negative real axis. We then redefine
$$F_{\tau, j}(\Delta) := F_j(\tau+\Delta), \qquad \tau \in \cz\setminus
\{0, -1, -2, \ldots\}$$
or in other words
\begin{align*}
F_{\tau,0}(\Delta) :&= f(\tau+\Delta), \\
F_{\tau,j}(\Delta) :&= 
F_{\tau,j-1}(\Delta+1) - (\tau+\Delta) F_{\tau,j-1}(\Delta), 
\qquad j=1, 2, \ldots \end{align*}
Using \eqref{recgam} we find
$$F_{\tau,j}(\Delta) = \sum_{i=1}^n \alpha_i\, \phi_i^j\;
\Gamma(\tau+\Delta+\phi_i), \qquad j=0, 1, 2, \ldots$$
If $\tau=0$ then $F_{0,j} = F_j(\Delta)$.
As soon as the samples at $\tau+\Delta+j$ are all
well-defined, we can start the algorithm for the
computation of the unknown linear
parameters $\alpha_i$ and the nonlinear parameters $\phi_i$,
We further introduce
$${^{\tau,k}_{\phantom{\tau,}1}}{\cal H}_n := \begin{pmatrix} 
F_{\tau,k} & \cdots & F_{\tau,k+n-1} \\ \vdots & \udots & \vdots \\
F_{\tau,k+n-1} & \cdots & F_{\tau,k+2n-2}
\end{pmatrix} .$$
\begin{theorem}
The matrix ${^{\tau,k}_{\phantom{\tau,}1}}{\cal H}_n$ is factored as
\begin{align*}
{^{\tau,k}_{\phantom{\tau,}1}}{\cal H}_n &= {\cal P}_n P_n Z_n {\cal P}_n^T, \\ 
{\cal P}_n &= \begin{pmatrix}
1 & \cdots  & 1 \\
\phi_1 & \cdots & \phi_n \\ \vdots & & \vdots \\
\phi_1^{n-1} & \cdots & \phi_n^{n-1} \end{pmatrix}, \\ 
Z_n &= \text{\rm diag} \left(\alpha_1\Gamma(\tau+\Delta+\phi_1),
\ldots, \alpha_n\Gamma(\tau+\Delta+\phi_n) \right), \\
P_n &= \text{\rm diag}(\phi_1^k, \ldots, \phi_n^k).
\end{align*}
\end{theorem}

\noindent{\it Proof. } With the matrix factorization given, the proof consists of an easy
verification of the matrix product with the matrix
${^{\tau,k}_{\phantom{\tau,}1}}{\cal H}_n$. \hspace{3pt}$\Box$ 
\smallskip

Filling the matrices ${^{\tau,0}_1}{\cal H}_n$ and ${^{\tau,1}_1}{\cal
H}_n$ requires the evaluation of $f(z)$
at $z=\tau+\Delta +j, j=0, \ldots, 2n-1$ which are points on a straight
line parallel with the real axis in the complex plane.
\smallskip

The nonlinear parameters $\phi_i$ are now obtained as the generalized
eigenvalues of
\begin{equation}\left( {^{\tau,1}_{\phantom{\tau,}1}}{\cal H}_n \right) v_i = \phi_i \left(
{^{\tau,0}_{\phantom{\tau,}1}}{\cal H}_n \right) v_i, \qquad i=1, \ldots, n,
\label{gep_G}\end{equation}
where the $v_i, i=1, \ldots, n$ are the right generalized eigenvectors.
Afterwards the linear parameters $\alpha_i$ are obtained from the
linear system of interpolation conditions 
$$\sum_{i=1}^n \left(\alpha_i \Gamma(\tau+\Delta+\phi_i) \right) \phi_i^j 
= F_{\tau,j}(\Delta), \qquad j=\tau, \ldots, \tau+2n-1,$$
by computing the coefficients $\alpha_i\Gamma(\tau+\Delta+\phi_i)$ and
dividing those by the function values $\Gamma(\tau+\Delta+\phi_i)$ which are
known because $\Delta,\tau$ and the $\phi_i, i=1, \ldots, n$ are known.
\smallskip

From Theorem 5 we find that for the generalized eigenvectors of
\eqref{gep_G} holds that ${^{\tau,0}_{\phantom{\tau,}1}}{\cal H}_n v_i$
is a scalar multiple of
$$\alpha_i \left( 1,
\phi_i, \ldots, \phi_i^{n-1} \right)^T.$$
This allows to validate the computation of the $\phi_i, i=1, \ldots, n$
obtained as generalized eigenvalues, if desired.

\subsection{Pochhammer basis connection}

Results on sparse polynomial interpolation using the Pochhammer basis
$(t)_m$
where usually the interpolation points are positive integers and $t \in \rz^+$,
were published in \cite{La.Sa:spa:95,Ka.Le:ear:03}, but no matrix pencil method
for its solution was presented. This can now easily be obtained using a
similar approach as for the gamma function. 
We consider more generally the interpolation of
$$f(z) = \sum_{i=1}^n \alpha_i \; (z)_{m_i}, \qquad z\in \cz\setminus
\{0\}, \quad m_i \in\nz.$$
For complex values $z$, the Pochhammer basis or rising factorial $(z)_m$ 
satisfies the recurrence relation
$$z\; \left[ (z+1)_m - (z)_m \right] = m\; (z)_m \;.$$
For real $\Delta$, a complex shift $\tau$ could shift the problem statement
away from  the negative real axis, as with the gamma function,
but it is much simpler here to immediately
consider $\Delta \in \cz \setminus \{0, -1, -2, \ldots\}$. Let
\begin{align*}
F_0 &:= F_0(\Delta) = f(\Delta), \\
F_j &:= F_j(\Delta) = 
\Delta\; \left[ F_{j-1}(\Delta+1) - F_{j-1}(\Delta) \right] = \sum_{i=1}^n
\alpha_i m_i^j\; (\Delta)_{m_i}, 
\qquad j = 1, 2, \ldots
\end{align*}
With the evaluations $F_j$ we fill the Hankel matrix
$${^k_1}H_n = \begin{pmatrix} F_k & F_{k+1} & \cdots & F_{k+n-1} \\
F_{k+1} & & & \\
\vdots & \udots & & \vdots \\
F_{k+n-1} & \cdots & & F_{k+2n-1}
\end{pmatrix} \;.$$
This Hankel matrix decomposes as in Theorem 5, but now with
\begin{align*}
{\cal P}_n &= \begin{pmatrix} 
1 & \cdots  & 1 \\
m_1 & \cdots & m_n \\ \vdots & & \vdots \\
m_1^{n-1} & \cdots & m_n^{n-1} \end{pmatrix}, \\ 
Z_n &= \text{\rm diag} \left(\alpha_1\; (\Delta)_{m_1},
\ldots, \alpha_n\; (\Delta)_{m_n} \right), \\
P_n &= \text{\rm diag}(m_1^k, \ldots, m_n^k).
\end{align*}
So the nonlinear parameters $m_i, i=1, \ldots, n$ are obtained as the
generalized eigenvalues of 
$${^1_1}H_n v_i = m_i\; {^0_1}H_n v_i, \qquad i=1, \ldots, n$$
where the $v_i$ are the right generalized eigenvectors, for which holds
that ${^0_1}H_n v_i$ is a multiple of the vector
$$\alpha_i\; (\Delta)_{m_i} \left( 1, m_i, \ldots, m_i^{n-1} \right)^T.$$
From the latter the estimates of the $m_i$ can be validated by computing
the quotient of successive entries in the vector.
The linear parameters $\alpha_i, i=1, \ldots, n$ are obtained from the
linear system
$$\sum_{i=1}^n \alpha_i m_i^j\; (\Delta)_{m_i} = F_j, \qquad j=0, \ldots,
2n-1.$$

\section{Numerical illustrations}

We present some examples to illustrate the main novelties of the
paper, including the multiscale facilities:
\begin{itemize}
\item an illustration of sparse interpolation 
by 
Gaussian distributions with fixed width but unknown peak locations; 
\item an illustration of the new generalized eigenvalue formulation for use
with several trigonometric functions and the sinc; 
\item an illustration of the use of the scale and shift strategy for the
supersparse interpolation of polynomials.
\end{itemize} 
As stated earlier, our focus is on the mathematical generalizations and 
not on the numerical issues. 

\subsection{Fixed width sparse Gaussian fitting}

Consider the expression
$$f(t) = \exp(-(t-5)^2) + 0.01 \exp(-(t-4.99)^2),$$
illustrated in Figure 1, with the parameters $\alpha_i, \phi_i \in \rz$.
From the plot it is not obvious that the signal has two peaks.
\smallskip

\vbox{
\centerline{\includegraphics[width=7truecm]{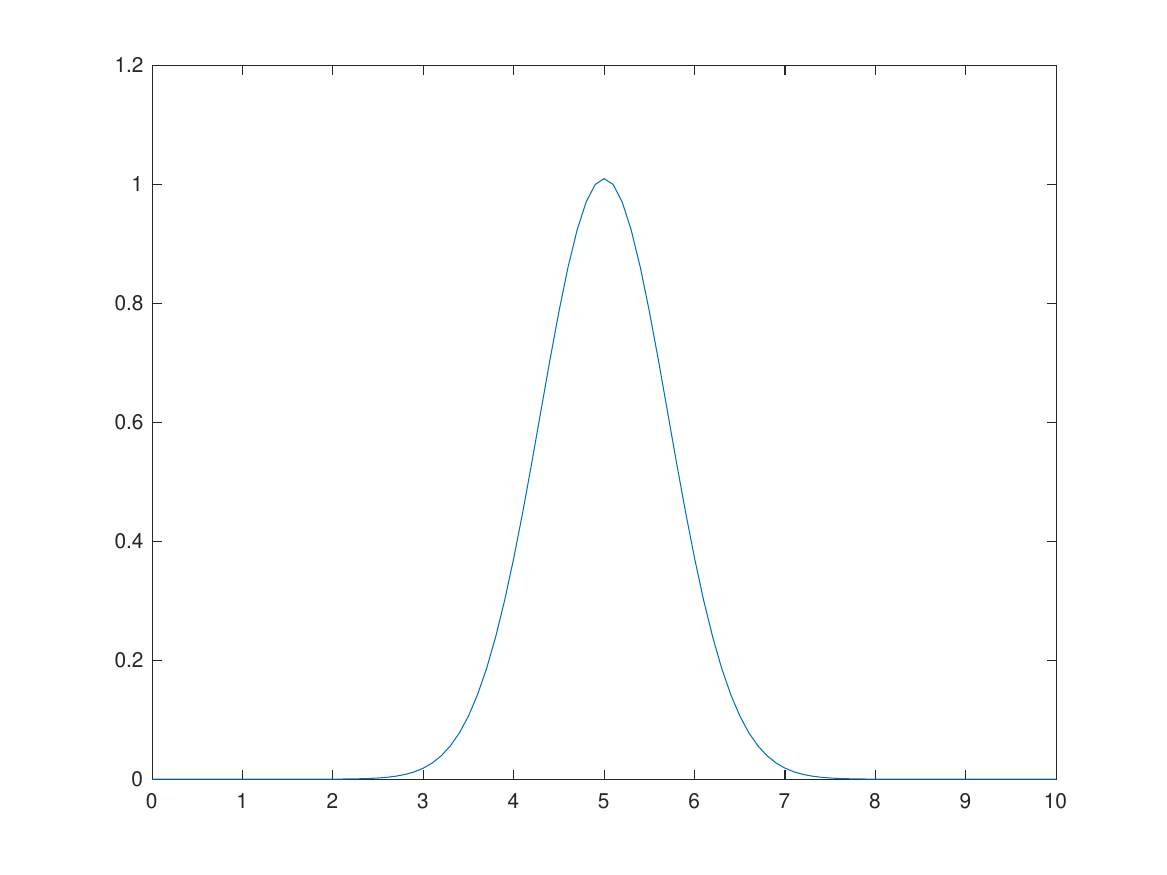}}
\smallskip
\centerline{\it Figure 1: Plot of Gaussian example function}
}
\smallskip

We first show the output of the widely used \cite{peakfituse}
Matlab state-of-the-art peak fitting program {\tt peakfit.m},
which calls an unconstrained nonlinear
optimization algorithm to decompose an overlapping peak signal into its
components \cite{peakfitdotm,peakfitbook}. 
\smallskip

In {\tt peakfit.m} the user needs to supply a guess for the number of
peaks and supply this as input. If one does not have any idea on the number of
peaks, the usual practice is to try different possibilities and compare
the corresponding results. Of course, a good estimate of the number of peaks
may lead to a good fit of the data. In addition to the peak position
$\phi_i$, its height $\alpha_i$ and width $w$, the program also returns a
goodness-of-fit ({\tt GOF}). 
\smallskip

The {\tt peakfit.m} algorithm can work without
assuming a fixed width, or the width can be passed as an argument. We do
the latter as our algorithm also assumes a known fixed peak width $w$.
\smallskip

Let $\Delta=0.1$ and let us collect 20 samples $f_0, f_1, \ldots, f_{19}$.
When passing the width to {\tt peakfit.m} and guessing the number of
peaks, then it returns for 1 peak the estimates
$$\phi_1=4.9998944950\ldots \qquad \alpha_1=1.0099771180\ldots \qquad
\text{\tt GOF} \approx 1.5 \times 10^{-5}.$$
For 2 peaks it returns
$$\begin{aligned} \phi_1 &= 4.9998944843\ldots \\ 
\phi_2 &= -1.5242538810\ldots 
\end{aligned} \qquad
\begin{aligned} 
\alpha_1 &= 1.0099776250\ldots \\
\alpha_2 &= 1.05 \times 10^{-13}
\end{aligned} \qquad 
\text{\tt GOF} \approx 5.2 \times 10^{-7}.$$
Since the result is still not matching our benchmark input parameters, 
let us push
further and supply 100 samples. Then for 1 peak {\tt peakfit.m} returns
$$\phi_1 = 4.9999009752\ldots \qquad \alpha_1 = 1.0099995049\ldots \qquad
\text{\tt GOF} \approx 2.5 \times 10^{-5}$$
and for 2 peaks we get
$$\begin{aligned} \phi_1 &= 4.9999945211\ldots \\ 
\phi_2 &= 4.9894206737\ldots 
\end{aligned} \qquad
\begin{aligned} 
\alpha_1 &= 1.0010660101\ldots \\
\alpha_2 &= 0.0089339897\ldots
\end{aligned} \qquad 
\text{\tt GOF} \approx 8.4 \times 10^{-9}.$$
From the latter 
experiment it is easy to formulate some desired features for a
new algorithm:
\begin{itemize}
\item built-in guess of the number of peaks in the signal,
\item and reliable output from a smaller number of samples.
\end{itemize} 
So let us investigate the technique developed in Section 5. Take
$\sigma=1$ and $\tau=0$ since there is no periodic component in the
Gaussian signal, which has only real parameters. With the 20 samples $f_0,
f_1, \ldots, f_{19}$ we define the samples
$F_j = \exp(j^2\Delta^2) f_j$ and compose the Hankel matrix
${_1^0}G_{10}$. Its singular value decomposition, illustrated in Figure 2, 
clearly
reveals that the rank of the matrix is 2 and so we deduce that there are
$n=2$ peaks. 
\smallskip

\vbox{
\centerline{\includegraphics[width=7truecm]{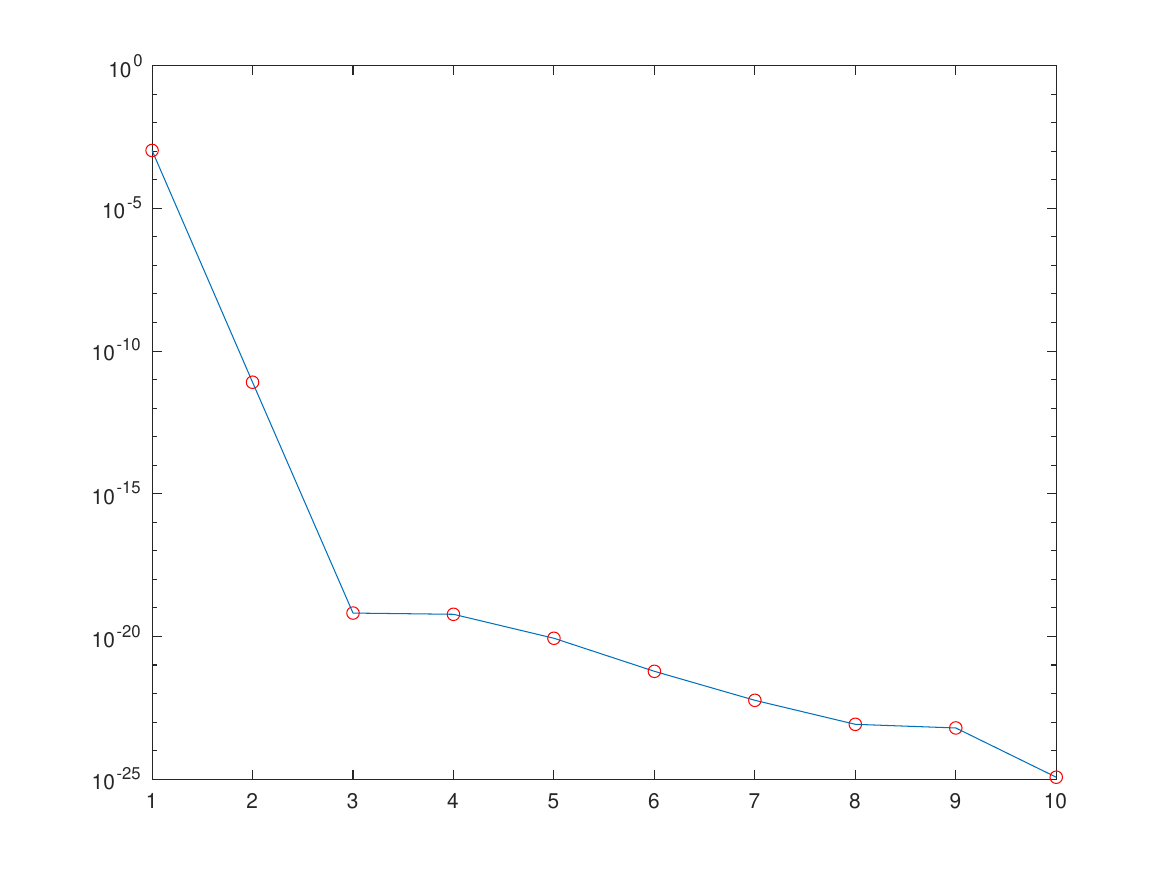}}
\smallskip
\centerline{\it Figure 2: Singular value log-plot of the matrix ${_1^0}G_{10}$}
}
\smallskip

From the 4 samples $F_0, F_1, F_2, F_3$ we obtain through Theorem 7
$$\begin{aligned} \phi_1 &= 4.9999999737\ldots \\ \alpha_1 &= 1.0000049866\ldots
\end{aligned} \qquad
\begin{aligned} \phi_2 &= 4.9899976207\ldots \\ \alpha_2 &= 0.0099950129\ldots
\end{aligned}$$ 
The new method clearly provides both an automatic guess of the number of
peaks and a reliable estimate of the signal parameters, all from only 20
samples. What remains to be done is to investigate the numerical behaviour
of the method on a large collection of different input signals, which
falls out of the scope of this paper where we provide the mathematical
details. 

\subsection{Sparse sinc interpolation}

Consider the function
$$f(t) = - 10 \text{sinc}(145.5t) + 20 \text{sinc}(149t) +
4 \text{sinc}(147.3t),$$
plotted in Figure 3,
which we sample at $t_j=j\pi/300$ for $j=0, \ldots, 19$. The singular
value decomposition of ${_1^0}B_{10}$ filled with the values
$t_j f_j$, of which the log-plot is shown in Figure 4 (left), 
reveals that $f(t)$ consists of 3 terms. Remember that the sparse sinc
interpolation problem with linear coefficients $\alpha_i$ and nonlinear
parameters $\phi_i$ transforms into a sparse sine interpolation problem
with linear coefficients $\alpha_i/\phi_i$ and samples $j\Delta f_j$.
\medskip

\vbox{
\centerline{\includegraphics[width=5truecm]{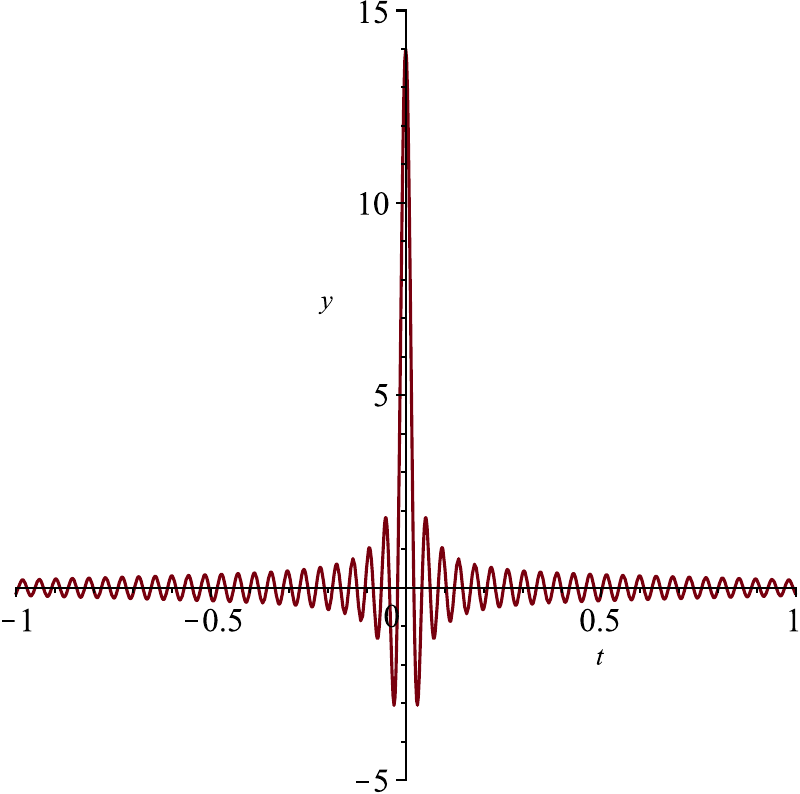}}
\smallskip
\centerline{\it Figure 3: The 3-term sparse sinc expression $f(t)$} 
}
\smallskip

The condition numbers of the matrices ${^0_1}B_3$ and ${^1_1}B_3$
appearing in the generalized eigenvalue problem
$$\left( {_1^1}B_3 \right) v_i = \cos(\phi_i\Delta) \left( {^0_1}B_3
\right), \qquad i=1, 2, 3$$
equal respectively $1.6 \times 10^7$ and $7.5 \times 10^6$.
To improve the conditioning of the structured matrix we choose
$\sigma=30, \tau=1$ and resample $f(t)$ at $t_j=30 j \pi/300 = j\pi/10$ 
for $j=0, \ldots, 5$. 
The singular values of ${^0_\sigma}B_{10}$ are graphed in Figure 4 (right)
and the condition numbers of ${^0_\sigma}B_3$ and ${^1_\sigma}B_3$ improve to
$1.1\times 10^3$ and $9.7\times 10^2$ respectively. 
\bigskip

\vbox{
\centerline{\hfill \includegraphics[width=4truecm]{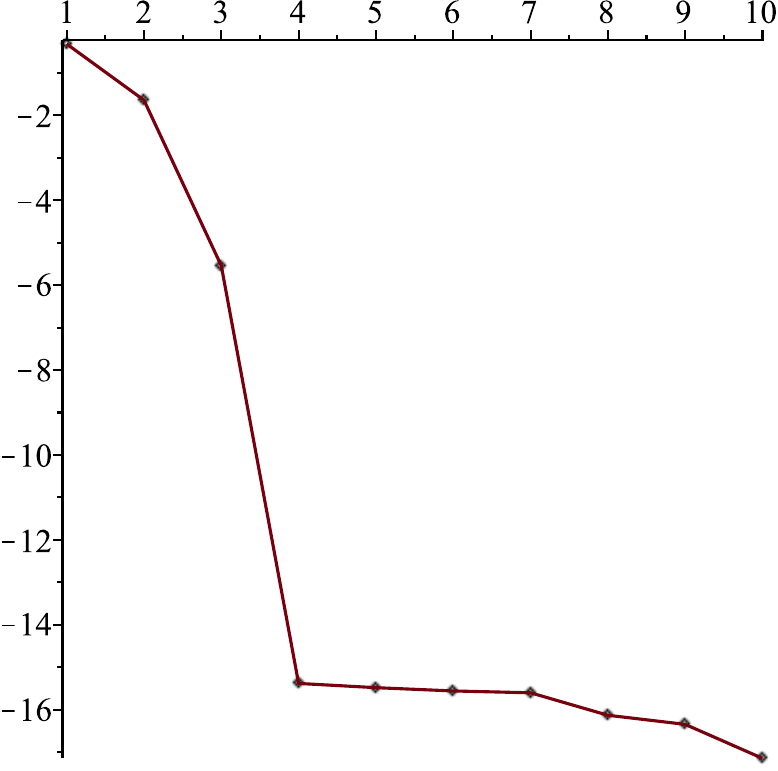} \hfill
\includegraphics[width=4truecm]{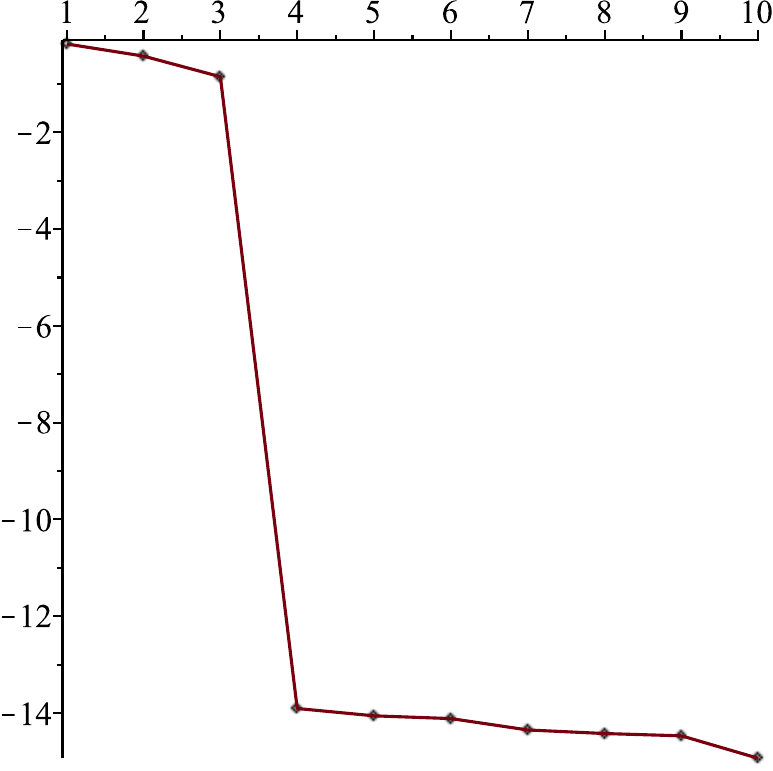} \hfill}
\smallskip

\centerline{\it Figure 4: 
Singular value log-plot of ${^0_1}B_{10}$ (left) and
${^0_{30}}B_{10}$ (right)}
}
\smallskip

The generalized eigenvalues of the matrix pencil
${^1_\sigma}B_3 -\lambda {^0_\sigma}B_3$ are given by
\begin{align*}
\cos(30\phi_1\Delta) &= -0.1564344650400536, \\
\cos(30\phi_2\Delta) &= -0.9510565162957546, \\
\cos(30\phi_3\Delta) &= -0.6613118653271576
\end{align*}
and with these we fill the matrix $W_3$ from Theorem 2. We solve
\eqref{sincA} for the
values $\alpha_i\sin(\phi_i\sigma\Delta)/\phi_i, i=1, 2, 3$ and further
compute for $j=1, \ldots, n,$
$${\alpha_i \over \phi_i} \sin(\phi_i j \sigma\Delta) = {\alpha_i \over
\phi_i} \sin(\phi_i (j-1)\sigma\Delta) \cos(\phi_i\sigma\Delta) +
\cos(\phi_i(j-1)\sigma\Delta) {\alpha_i \over \phi_i} 
\sin(\phi_i\sigma\Delta).$$ 
At this point the matrix $U_3$ from Theorem 2 can be filled and the
$\cos(\phi_i\tau\Delta)$ can be computed from \eqref{sincB}, with the
right hand side filled with 
the additional samples $F(31\Delta), F(61\Delta),
F(91\Delta)$, where
$$F_{\tau+j\sigma} = {\Delta \over 2} (\tau+j\sigma) f_{\tau+j\sigma} +
{\Delta \over 2} (-\tau+j\sigma) f_{-\tau+j\sigma}.$$
Since $\tau=1$ we obtain the $\phi_i$ directly from the values
$\cos(\phi_i\tau\Delta)$:
$\phi_1 = 145.5000000000$, 
$\phi_2 = 149.0000000000$,  
$\phi_3 = 147.3000000000$.
The linear coefficients $\alpha_i$ are given by
$$\alpha_i = \phi_i {\alpha_i \sin(\phi_i\sigma\Delta)/\phi_i \over
\sin(\phi_i\sigma\Delta)},$$
resulting in
$\alpha_1 = -9.999999999991,   
\alpha_2 = 19.99999999978,   
\alpha_3 = 4.000000000089.$

\subsection{Supersparse Chebyshev interpolation}

We consider the polynomial  
$$f(t) = 2 T_6(t) + T_7(t) + T_{39999}(t),$$
which is clearly supersparse when expressed in the Chebyshev basis.
We sample $f(t)$ at $t_j=\cos(j\Delta)$ where $\Delta=\pi/100000$ with
$M=50000$.
The first challenge is maybe to retrieve an indication of the sparsity $n$. 
\smallskip

Take $\sigma=1$ and collect 15 samples $f_j, j=0, \ldots, 14$ to form the
matrix ${^0_1}C_8$. From its singular value decomposition, computed in
double precision arithmetic and illustrated on the log-plot in Figure 5, 
one may erroneously conclude that $f(t)$ has only 2 terms, a consequence of the
fact that, relatively speaking, the degrees $m_1=6$ and $m_2=7$ are close to
one another and appear as one cluster.
\medskip

\vbox{
\centerline{\includegraphics[width=4truecm]{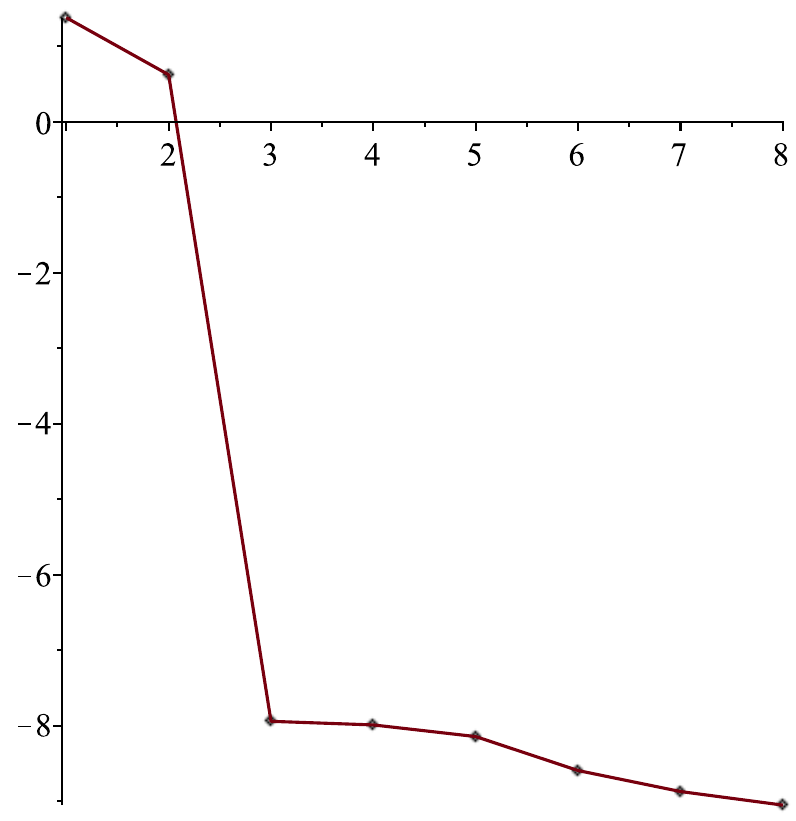}}
\smallskip

\centerline{\it Figure 5: Log-plot of singular values of ${^0_1}C_8$}
}
\smallskip

Imposing a 3-term model to $f(t)$ instead of the erroneously suggested
2-term one, does not improve the computation as the
matrix ${^0_1}C_8$ is ill-conditioned with a condition number of the order
of $10^{10}$. So 3 generalized
eigenvalues cannot be extracted reliably from the samples. For
completeness we mention the unreliable double precision 
results, rounded to integer values: $m_1 = 6, m_2 = 39999, m_3 = 25119$.
\smallskip

Now choose $\sigma=3125$ and $\tau=16$. The singular value decomposition
of ${^0_{3125}}C_8$, shown on the log-plot in Figure 6,
reveals that $f(t)$ indeed consists of 3 terms. Also, the conditioning of
the involved matrix ${^0_{3125}}C_8$ has improved to the order of $10^3$.
\medskip

\vbox{
\centerline{\includegraphics[width=4truecm]{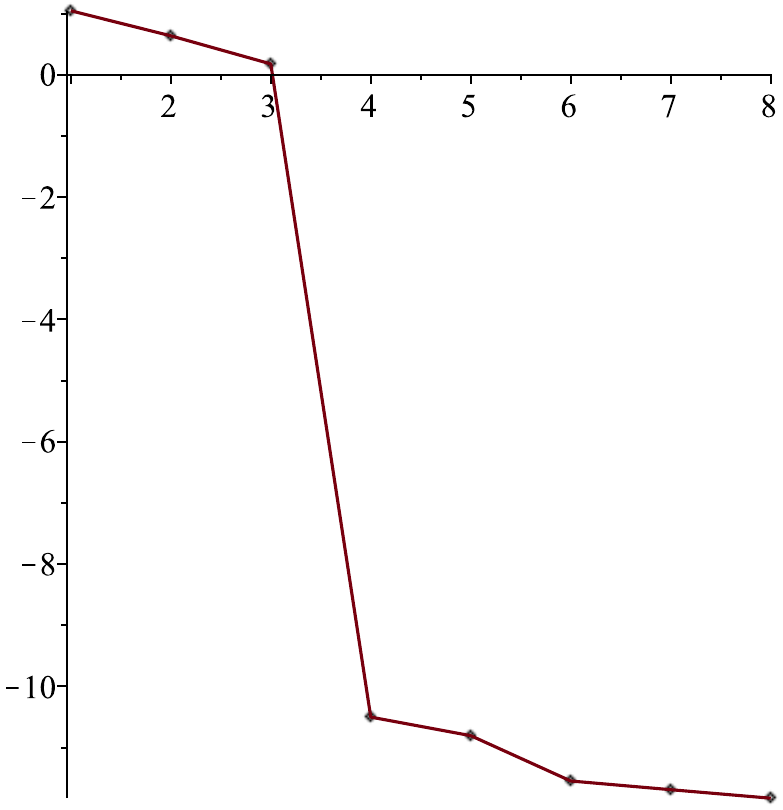}}
\smallskip

\centerline{\it Figure 6: Log-plot of singular values of ${^0_{3125}}C_8$}
}
\smallskip

The distinct
generalized eigenvalues extracted from the matrix pencil ${^1_{3125}}C_3
-\lambda {^0_{3125}}C_3$ are given by
\begin{align*}
\cos(m_1\sigma\Delta) &= 0.9999999204093383, \\ 
\cos(m_2\sigma\Delta) &= -0.8089800617792506, \\ 
\cos(m_3\sigma\Delta) &= -0.007490918959382487. 
\end{align*}
From 3 shifted samples at the arguments $t_{\tau+j\sigma}, j=0, 1, 2$
we obtain
\begin{align*}
\cos(m_1\tau\Delta) &= -0.8084256802389809, \\ 
\cos(m_2\tau\Delta) &= 0.9999752362021560, \\ 
\cos(m_3\tau\Delta) &= 0.9999818099296417. 
\end{align*}
Building the sets $S_i$ and $T_i$ for $i=1, 2, 3$ as indicated in Section
4.1, and rounding the result to the nearest integer, does unfortunately not
provide singletons for $S_1\cap T_1, S_2\cap T_2, S_3\cap T_3$. We
consequently need to consider a second shift, for which we choose
$\sigma+\tau=3141$. With this choice
we only need to add the evaluation of $f(\tau+n\sigma)$
to proceed and compute
\begin{align*}
\cos(m_1(\sigma+\tau)\Delta) &= -0.6780621808989576, \\
\cos(m_2(\sigma+\tau)\Delta) &= 0.1881836009619241, \\
\cos(m_3(\sigma+\tau)\Delta) &= 0.3771037932233129.
\end{align*}
Finally, intersecting each $S_i\cap T_i, i=1, 2,3$
with the solutions provided by the second shift, delivers the correct 
$m_1 =6, 
m_2 = 7, 
m_3 = 39999.$

\section{Conclusion}

Let us summarize the sparse interpolation formulas 
obtained in the preceding sections in a table.
For each parameterized
univariate function $g(\phi_i; t)$ we list in the columns 1 to 4:
\begin{enumerate}
\item the minimal number of
samples required to solve the sparse interpolation without running
into ambiguity problems, meaning for the choice $\sigma=1$,

\item the minimal number of samples required for the choice $\sigma >1$ (if
applicable), 
thereby involving a shift $\tau \not= 0$ to restore uniqueness
of the solution,

\item the linear matrix pencil 
$(A, B)$ in the generalized eigenvalue formulation
$Av_i = \lambda_i Bv_i$ of the
sparse interpolation problem involving the $g(\phi_i;t)$,  

\item the generalized eigenvalues in terms of $\tau$, as they can be read
directly from the structured matrix factorizations presented in the
theorems 1--5,

\item the information that can be computed from the associated
generalized eigenvectors, as indicated at the end of each (sub)section.
\end{enumerate} 

\begin{center}
\begin{tabular}{ ? c | c | c ? c | c | c ?}
\Xhline{3\arrayrulewidth}
 $g(\phi_i; t)$ &  \multicolumn{2}{c?}{\# samples} &  pencil$^*$ $(A, B)$ & $\lambda_i$ & $Bv_i$ \\ 
 \hline
  & $\sigma=1$ & $\sigma>1$ &  \multicolumn{3}{c?}{ } \\  
  \hline
 $\exp(\phi_it)$ & $2n$ & $3n$ & $\left( {_\sigma^\tau}H_n, {_\sigma^0}H_n \right)$ & $\exp(\phi_i\tau\Delta)$ & $\alpha_i, \exp(\phi_i\sigma\Delta)$ \\  
 \hline
 $\cos(\phi_i t)$ & $2n$ & $4n$  & $\left( {_\sigma^\tau}C_n, {_\sigma^0}C_n \right)$ & $\cos(\phi_i\tau\Delta)$ &  $\alpha_i, \cos(\phi_i\sigma\Delta)$\\
 \hline
 $\sin(\phi_i t)$ & $2n$ & $4n+2$  & $\left( {_\sigma^\tau}B_n, {_\sigma^0}B_n \right)$ & $\cos(\phi_i\tau\Delta)$ & $\alpha_i, \sin(\phi_i\sigma\Delta)$ \\
 \hline
 $\cosh(\phi_i t)$ & $2n$ & $4n$  & $\left( {_\sigma^\tau}C_n^*, {_\sigma^0}C_n^* \right)$ & $\cosh(\phi_i\tau\Delta)$  & $\alpha_i, \cosh(\phi_i\sigma\Delta)$  \\
 \hline
 $\sinh(\phi_i t)$ & $2n$ & $4n+2$  & $\left( {_\sigma^\tau}B_n^*, {_\sigma^0}B_n^* \right)$ & $\cosh(\phi_i\tau\Delta)$  & $\alpha_i, \sinh(\phi_i\sigma\Delta)$ \\
 \hline
 $T_{m_i}(t)$ & $2n$ & $4n$  &  $\left( {_\sigma^\tau}C_n, {_\sigma^0}C_n \right)$ & $T_{m_i}(\cos \tau\Delta)$ &  $\alpha_i, T_{m_i}(\cos \sigma\Delta)$\\
 \hline
  $S_{m_i}(t)$ & $2n+1$ & $4n+2$  & $\left( {_\sigma^\tau}K_n, {_\sigma}J_n \right)$  & $S_{m_i}(\sin^2 \tau\Delta)$ & $\alpha_i, S_{m_i}(\sin^2\sigma\Delta)$ \\
  \hline
  $\text{sinc}(\phi_i t)$ & $2n$ & $4n+2$  & $\left( {_\sigma^\tau}B_n, {_\sigma^0}B_n \right)$ & $\cos(\phi_i\tau\Delta)$ & $\alpha_i, \sin(\phi_i\sigma\Delta)$ \\
  \hline
  $\Gamma(z+\phi_i)$ & $2n$ & {\tt X}  & $\left( {_\sigma^{\tau,1}}{\cal H}_n, {_\sigma^{\tau,0}}{\cal H}_n \right)$ & $\phi_i$  & $\alpha_i, \phi_i$  \\
  \hline
  $\exp(-(t-\phi_i)^2)$ & $2n$ & $3n$  & $\left( {_\sigma^\tau}G_n, {_\sigma^0}G_n \right)$ & $\exp(2\phi_i\tau\Delta)$ & $\alpha_i, \exp(2\phi_i\sigma\Delta)$  \\
 \Xhline{3\arrayrulewidth}
\end{tabular}
\smallskip

\thefootnote{$*$ with $\cos(\cdot)$ replaced by $\cosh(\cdot)$ and
$\sin(\cdot)$ replaced by $\sinh(\cdot)$ in the hyperbolic case}
\end{center}

\section*{Appendix}

To reconstruct a function of the form
$$f(t) = \sum_{i=1}^n \alpha_i \cos(\phi_i t)$$
from equidistantly collected samples $f_j$ at $t=j\Delta$, in other words to
recover the unknown parameters $\phi_i$ and coefficients $\alpha_i$,
the sampling step $\Delta$ needs to satisfy the Shannon-Nyquist constraint
$$\Delta < \pi/ \max_{i=1, \ldots, n} |\phi_i|.$$
Since we do not distinguish $\phi_i$ from $-\phi_i$ in this case, we can
simply drop the sign information in $\phi_i$ from here on 
and write $\Delta = \pi/R$ with
$$0 \le \max_{i=1, \ldots, n} \phi_i < R.$$
The challenge we consider now is to retrieve the parameters $\phi_i$ and
coefficients $\alpha_i$ from sub-Nyquist rate collected samples $f_{j\sigma}$ at
$t=j\sigma\Delta$ with $\sigma>1$ and the shifted evaluations
$f_{j\sigma +\tau}$ at $t=(j\sigma+\tau)\Delta$ with
$\gcd(\sigma,\tau)=1$. In Section 3.1 we describe how for $i=1,
\ldots, n$ the values
\begin{align*}
C_{i,\sigma} :=& \cos(\phi_i\sigma\Delta), \\
C_{i,\tau} :=& \cos(\phi_i\tau\Delta) 
\end{align*}
are obtained. The aim is to extract the correct value for $\phi_i$ from the
knowledge of the evaluations $C_{i,\sigma}$ and $C_{i,\tau}$, particularly
when $(\sigma \Delta) \max_{i=1, \ldots,
n}\phi_i \ge \pi$ and the parameter $\phi_i$ cannot be obtained uniquely
from $C_{i,\sigma}$ alone. 
We now discuss the unique identification of this parameter
$\phi_i$ and in doing so we further drop the index $i$. Let us denote
\begin{equation}
\label{arccos}
\begin{aligned}
A_{\sigma} :=& \Arccos (C_{\sigma}), \\
A_{\tau} :=& \Arccos (C_{\tau})
\end{aligned} \end{equation}
where $\Arccos(\cdot) \in [0,\pi]$ 
indicates the principal value of the inverse cosine
function. Knowing that $0 \le A_\sigma, A_\tau \le \pi$ and that 
$0 \le \phi\sigma\Delta < \sigma\pi$, we find that all possible 
positive arguments $\phi\sigma\Delta$ of $C_\sigma$ are in ${\cal
A}_{\sigma,1} \cup {\cal A}_{\sigma,2}$ with
\begin{align*}
{\cal A}_{\sigma,1} :=& \{A_\sigma + 2\pi\ell \mid 0 \le \ell \le
\lceil\sigma/2\rceil -1 \}, \\ 
{\cal A}_{\sigma,2} :=& \{(2\pi-A_\sigma)\sgn(A_\sigma) + 2\pi\ell \mid 
0 \le \ell \le \lceil\sigma/2\rceil -1 \},
\end{align*}
where $\sgn(a_\sigma)=+1$ for $0<A_\sigma \le \pi$ and $\sgn(0)=0$.
The set ${\cal A}_{\sigma,1} \cup {\cal A}_{\sigma,2}$ may even contain 
some candidate arguments of $C_\sigma$ that do not satisfy the bounds, 
but this does not create a problem
in the identification of the correct $\phi<R$. Along the same lines, sets
${\cal A}_{\tau,1}$ and ${\cal A}_{\tau,2}$ can be constructed.
\smallskip

We further denote
\begin{equation}
\label{mulcos}
\begin{aligned}
\phi_{\sigma} :=& A_\sigma/(\sigma\Delta) = {A_\sigma R \over \sigma \pi}, \\
\phi_{\tau} :=& A_\tau/(\tau\Delta) = {A_\tau R \over \tau \pi}.
\end{aligned}
\end{equation} 
Then the possible solutions for $\phi$ 
to $C_\sigma=\cos(\phi\sigma\Delta)$ are in
$\Phi_{\sigma,1} \cup \Phi_{\sigma,2}$ where
\begin{align}
\Phi_{\sigma,1} :=& \{\phi_\sigma + 2R \ell/\sigma \mid 0 \le \ell \le
\lceil\sigma/2\rceil -1 \} \cap [0,R), \label{S1} \\ 
\Phi_{\sigma,2} :=& \{(2R/\sigma-\phi_\sigma)\sgn(\phi_\sigma) + 2R \ell/\sigma \mid 0 \le \ell \le
\lceil\sigma/2\rceil -1 \} \cap [0,R). \label{S2}
\end{align}
Analogously, the possible solutions to $C_\tau=\cos(\phi\tau\Delta)$ are in 
$\Phi_{\tau,1} \cup \Phi_{\tau,2}$ where
\begin{align}
\Phi_{\tau,1} :=& \{\phi_\tau + 2R \ell/\tau \mid 0 \le \ell \le
\lceil\tau/2\rceil -1 \} \cap [0,R), \label{T1} \\
\Phi_{\tau,2} :=& \{(2R/\tau-\phi_\tau)\sgn(\phi_\tau) + 2R \ell/\tau \mid 
0 \le \ell \le \lceil\tau/2\rceil -1 \} \cap [0,R). \label{T2}
\end{align}
One statement is obvious: whatever the choice for $\sigma$ and $\tau$,
both $\Phi_{\sigma,1} \cup \Phi_{\sigma,2}$ and $\Phi_{\tau,1} \cup
\Phi_{\tau,2}$ contain the unknown value for $\phi$ which produced
$C_\sigma$ and $C_\tau$. What remains open is the question whether
$\left( \Phi_{\sigma,1}\cup \Phi_{\sigma,2} \right) \cap 
\left( \Phi_{\sigma,1}\cup \Phi_{\sigma,2} \right)$
is a singleton. And in case it is not, we want to find an
algorithm that can identify the correct $\phi$.
\smallskip

When either $\phi_\sigma=0$ or $\phi_\sigma = R/\sigma$ the 
sets $\Phi_{\sigma,1}$ and
$\Phi_{\sigma,2}$ coincide. And similarly for $\phi_\tau$. On the other
hand, if these sets do not coincide, they are disjoint.
So the true value for the unknown parameter $\phi$ 
can belong to any of the intersections $\Phi_{\sigma,1} \cap
\Phi_{\tau,1}, \Phi_{\sigma,1} \cap \Phi_{\tau,2}, \Phi_{\sigma,2} \cap 
\Phi_{\tau,1}, \Phi_{\sigma,2} \cap \Phi_{\tau,2}$. 
A sequence of lemmas
will lead to the conclusion that the four intersections do not deliver more
than two distinct elements. Thereafter we indicate how to identify the
only true value for the unkown $\phi$.
\smallskip

\begin{lemma}
$i,j \in \{1, 2\}: \Phi_{\sigma,i} \cap \Phi_{\tau,j} \not= \emptyset 
\Longrightarrow \#\left(
\Phi_{\sigma,i} \cap \Phi_{\tau,j} \right) = 1.$
\end{lemma}
\smallskip

\noindent{\it Proof. }
Without loss of generality we prove the statement for $i=1=j$, by
contraposition. The proof of the other cases is entirely similar. From
$\Phi_{\sigma,1} \cap \Phi_{\tau,1} \not= \emptyset$ and containing at least
two elements, we then find that
\begin{multline*}\exists 0 \le \ell_1, \ell_2 \le \lceil \sigma/2 \rceil -1, 0 \le k_1, k_2 
\le \lceil \tau/2 \rceil -1, \ell_1\not=\ell_2, k_1\not=k_2: \\
\left\{ \begin{aligned} \phi_\sigma + \ell_1 2R/\sigma &= \phi_\tau + k_1 2R/\tau \\
\phi_\sigma + \ell_2 2R/\sigma &= \phi_\tau + k_2 2R/\tau.
\end{aligned} \right.\end{multline*}
This leads to
$${\ell_1 -\ell_2 \over k_1 - k_2} = {\sigma \over \tau}$$
which is a contradiction because $|\ell_1-\ell_2| < \sigma, |k_1-k_2|< \tau$
and $\gcd(\sigma,\tau)=1$. \hspace{3pt}$\Box$  
\smallskip

When the sets $\Phi_{\sigma,1}$ and $\Phi_{\sigma,2}$ coincide and the sets 
$\Phi_{\tau,1}$ and $\Phi_{\tau,2}$ do as well, then that unique
intersection is $\phi=\phi_\sigma=\phi_\tau=0$. Because
$\gcd(\sigma,\tau)=1$, other common elements coming from either
$\phi_\sigma=2R/\sigma$ or $\phi_\tau=2R/\tau$ cannot exist.
\smallskip

We now continue with 
the situation where either the sets in \eqref{S1} and \eqref{S2} or
the sets in \eqref{T1} and \eqref{T2} do not coincide, so that there are always
at least 3 distinct sets in  the running.
Without loss of
generality, we assume that a common element belongs to $\Phi_{\sigma,1} \cap
\Phi_{\tau,1}$ and we build our reasoning from there. 
\smallskip

\begin{lemma} $\Phi_{\sigma,1} \cap \Phi_{\tau,1} \not= \emptyset
\Longrightarrow \Phi_{\sigma,2} \cap \Phi_{\tau,2} = 
\emptyset$. 
\end{lemma}
\smallskip

\noindent{\it Proof. } We know that either $\Phi_{\sigma,1} \cap \Phi_{\sigma,2} \not= \emptyset$
or $\Phi_{\tau,1} \cap \Phi_{\tau,2} \not= \emptyset$ and possibly both,
so that $\Phi_{\sigma,1} \cap \Phi_{\tau,1} \not= \Phi_{\sigma,2} \cap 
\Phi_{\tau,2}$.
Again by contraposition, we suppose that $\Phi_{\sigma,2} \cap \Phi_{\tau,2}
\not= \emptyset$ and so
\begin{multline*}
\exists 0 \le \ell_1, \ell_2 \le \lceil \sigma/2 \rceil -1, 0 \le k_1, k_2
\le \lceil \tau/2 \rceil -1: \\ 
\left\{ \begin{aligned} \phi_1 &= \phi_\sigma + \ell_1 2R/\sigma = \phi_\tau + 
k_1 2R/\tau \in \Phi_{\sigma, 1}\cap\Phi_{\tau,1}, \\
\phi_2 &= 2R/\sigma-\phi_\sigma + \ell_2 2R/\sigma = 2R/\tau-\phi_\tau + 
k_2 2R/\tau \in \Phi_{\sigma, 2}\cap\Phi_{\tau,2}.
\end{aligned} \right. \end{multline*}
From this we obtain
$${1+\ell_1+\ell_2 \over 1+k_1+k_2} = {\sigma \over \tau}$$
which can only be true when $1+\ell_1+\ell_2 = \sigma$ and $1+k_1+k_2=\tau$.
Then
$$\phi_1+\phi_2 = (1+\ell_1+\ell_2)2R/\sigma = 2R,$$
which contradicts $0 \le \phi_1, \phi_2 <R$. \hspace{3pt}$\Box$ 
\smallskip

While, assuming $\Phi_{\sigma,1}\cap\Phi_{\tau,1} \not= \emptyset$, we have
seen in Lemma 1 that this intersection is a singleton, and we have seen
in Lemma 2  that then $\Phi_{\sigma,2}\cap\Phi_{\tau,2} = \emptyset$, we know
nothing so far about the other two intersections
$\Phi_{\sigma,1}\cap\Phi_{\tau,2}$ and $\Phi_{\sigma,2}\cap\Phi_{\tau,1}$. 
\smallskip

\begin{lemma} $\Phi_{\sigma,1}\cap\Phi_{\tau,1} \not= \emptyset
\Longrightarrow \neg \left( \Phi_{\sigma,2} \cap \Phi_{\tau,1}
\not= \emptyset \wedge
\Phi_{\sigma,1} \cap \Phi_{\tau,2} \not= \emptyset \right)$.
\end{lemma}
\smallskip

\noindent{\it Proof. } By contraposition we assume that 
\begin{multline*}
\exists 0 \le \ell_1, \ell_2 \le \lceil \sigma/2 \rceil -1, 0 \le k_1, k_2
\le \lceil \tau/2 \rceil -1: \\
\left\{ \begin{aligned} \phi_1 &= 2R/\sigma-\phi_\sigma + \ell_1 2R/\sigma = \phi_\tau +
k_1 2R/\tau \in \Phi_{\sigma, 2}\cap\Phi_{\tau_1}, \\
\phi_2 &= \phi_\sigma + \ell_2 2R/\sigma = 2R/\tau-\phi_\tau +
k_2 2R/\tau \in \Phi_{\sigma, 1}\cap\Phi_{\tau,2}.
\end{aligned} \right. \end{multline*}
This leads to 
$${1+\ell_1+\ell_2 \over 1+k_1+k_2} = {\sigma \over \tau},$$
which again implies $1+\ell_1+\ell_2=\sigma$ and $1+k_1+k_2=\tau$. Since
$$\phi_1+\phi_2 = (1+\ell_1+\ell_2)2R/\sigma = 2R,$$
this contradicts $0 \le \phi_1, \phi_2 <R$. \hspace{3pt}$\Box$ 
\smallskip

We have built our sequence of proofs from Lemma 2 on,
without loss of generality, on the fact that $\Phi_{\sigma, 1}
\cap \Phi_{\tau,1} \not= \emptyset$ and the fact that $\Phi_{\sigma,1}$ and
$\Phi_{\sigma,2}$ on the one hand and $\Phi_{\tau,1}$ and $\Phi_{\tau,2}$
on the other do not collide at the same time. 
Finally, from Lemma 3 we know that
($\Phi_{\sigma,1}\cap\Phi_{\tau,1}
\not= \emptyset$ and $\Phi_{\sigma,2}\cap\Phi_{\tau,1} \not= \emptyset$) or
($\Phi_{\sigma,1}\cap\Phi_{\tau,1} \not= \emptyset$ and 
$\Phi_{\sigma,1}\cap\Phi_{\tau,2} \not= \emptyset$)
cannot occur concurrently, but either one of these cases remains possible.
\smallskip

In general, when at least 3 of the 4 sets $\Phi_{\sigma, 1},
\Phi_{\sigma, 2}, \Phi_{\tau,1}, \Phi_{\tau,2}$ are distinct,
then at most 2 of the 4 intersections
$$\Phi_{\sigma,i} \cap \Phi_{\tau,j}, \qquad 1 \le i,j \le 2$$
are nonempty, with each of the nonempty intersections being a singleton.
Further down we illustrate the actual existence of a case, where
two intersections are nonempty and consequently the
true value of the unknown $\phi$ cannot be identified from the evaluations
$\cos(\phi\sigma\Delta)$ and $\cos(\phi\tau\Delta)$ with
$\gcd(\sigma,\tau)=1$. 
\smallskip

In this case we need to collect a third value
$C_\rho:=\cos(\phi\rho\Delta)$ with $\gcd(\sigma,\rho)=1$ and 
$\gcd(\tau,\rho)=1$.
With $A_\rho$ and $\phi_\rho$ defined as in \eqref{arccos} 
and \eqref{mulcos}, and $\Phi_{\rho,1}$ and $\Phi_{\rho,2}$ defined
as in \eqref{S1} and \eqref{S2}, we know, as before, that 
$\Phi_{\rho,1}\cup\Phi_{\rho,2}$ contains the correct value for $\phi$.
We also know, because of the remark formulated
after the proof of Lemma 1, that at least 5 of the 6 involved
sets ${\Phi}_{\sigma,1}, {\Phi}_{\sigma,2}, {\Phi}_{\tau,1}, 
{\Phi}_{\tau,2}, {\Phi}_{\rho,1}, {\Phi}_{\rho,2}$ are distinct unless
$\phi=0$.
\smallskip

We now inspect
\begin{equation}
\label{rho}
\begin{aligned}
\left[\cup_{i, j=1}^2 \left( \Phi_{\sigma, i} \cap \Phi_{\tau, j} \right) 
\right] 
&\cap \left( \Phi_{\rho,1} \cup \Phi_{\rho,2} \right) \\
&= \left( \cup_{k=1}^2 \Phi_{\sigma, i_1} \cap \Phi_{\tau,j_1} \cap
\Phi_{\rho,k} \right) \cup
\left( \cup_{k=1}^2 \Phi_{\sigma, i_2} \cap \Phi_{\tau,j_2} \cap
\Phi_{\rho,k} \right) 
\end{aligned} \end{equation}
where $i_1, j_1, i_2, j_2$ index the subsets that produce the nonempty
intersections of the relatively prime pair $\sigma$ and $\tau$,
with either $i_1 \not= i_2$ or $j_1 \not= j_2$ but not both. 
We have built our sequence of proofs, 
without loss of generality, on the fact that 
$i_1=1, j_1=1$ and have found that it is then possible that $i_2=2, j_2=1$. 
We now continue the proofs from that case and inspect the 4 new intersections
in \eqref{rho}.
\smallskip

\begin{lemma} 
$\Phi_{\rho,1} \cap \Phi_{\rho,2} =\emptyset \wedge
\Phi_{\sigma,1}\cap\Phi_{\tau,1}\cap\Phi_{\rho,1} \not=
\emptyset \Longrightarrow \Phi_{\sigma,1}\cap\Phi_{\tau,1}\cap\Phi_{\rho,2} =
\emptyset.$
\end{lemma}
\smallskip

\noindent{\it Proof. } From Lemma 1, we know that $\Phi_{\sigma, 1} \cap \Phi_{\tau,1}$ is a singleton. 
If that unique element also belongs to $\Phi_{\rho,1}$ then it cannot
belong to $\Phi_{\sigma,1}\cap\Phi_{\tau,1}\cap\Phi_{\rho,2}$ when
$\Phi_{\rho,1}$ and $\Phi_{\rho,2}$ are disjoint. \hspace{3pt}$\Box$ 
\smallskip

\begin{lemma} 
$\Phi_{\sigma,1} \cap \Phi_{\sigma,2} =\emptyset \wedge
\Phi_{\sigma,1}\cap\Phi_{\tau,1}\cap\Phi_{\rho,1} \not= \emptyset
\Longrightarrow
\Phi_{\sigma,2}\cap\Phi_{\tau,1}\cap\Phi_{\rho,1} = \emptyset$. 
\end{lemma}
\smallskip

\noindent{\it Proof. } 
From Lemma 1, we know that $\Phi_{\tau, 1} \cap \Phi_{\rho,1}$ is a singleton.
If that unique element also belongs to $\Phi_{\sigma,1}$ then it cannot
belong to $\Phi_{\sigma,2}\cap\Phi_{\tau,1}\cap\Phi_{\rho,1}$ when
$\Phi_{\sigma,1}$ and $\Phi_{\sigma,2}$ are disjoint. \hspace{3pt}$\Box$ 
\smallskip

As a consequence of the Lemmas 4 and 5, the unique true $\phi$ is identified
in 
$$\Phi_{\sigma,1}\cap\Phi_{\tau,1}\cap\Phi_{\rho,1} \text{ or } 
\Phi_{\sigma,2}\cap\Phi_{\tau,1}\cap\Phi_{\rho,2}.$$
\begin{lemma}
$\#\left[\left( \Phi_{\sigma,1}\cap\Phi_{\tau,1}\cap\Phi_{\rho,1} \right) \cup
\left( \Phi_{\sigma,2}\cap\Phi_{\tau,1}\cap\Phi_{\rho,2} \right) \right] = 1$.
\end{lemma}
\smallskip

\noindent{\it Proof. } 
We know that either $\phi=0$ is the unique element in the intersections
or at least 2 of the intersections 
$\Phi_{\sigma,1} \cap \Phi_{\sigma,2}, \Phi_{\tau,1} \cap \Phi_{\tau,2}$,
\hbox{$\Phi_{\rho,1} \cap \Phi_{\rho,2}$} are empty. So either $\Phi_{\sigma,1}
\cap \Phi_{\sigma,2} =\emptyset$ or $\Phi_{\rho,1} \cap \Phi_{\rho,2} =
\emptyset$. When applying Lemma 2 to the
pair $(\sigma,\rho)$ instead of $(\sigma,\tau)$ the set
$\Phi_{\sigma,2}\cap\Phi_{\rho,2}=\emptyset$ if the set
$\Phi_{\sigma,1}\cap\Phi_{\rho,1} \not= \emptyset$. Therefore two distinct
elements in respectively
$\Phi_{\sigma,1}\cap\Phi_{\rho,1}$ and $\Phi_{\sigma,2}\cap\Phi_{\rho,2}$
cannot coexist and solve $C_\tau=\cos(\phi\tau\Delta)$. \hspace{3pt}$\Box$ 
\smallskip

So the unknown parameter $\phi$ is identified uniquely from at most 3
values $C_\sigma, C_\tau, C_\rho$
with $\sigma, \tau, \rho$ all mutually prime. 
An easy choice for $\rho$ is
$\rho=\sigma+\tau$ as this minimizes the number of additional samples as
explained in Section 3, and also
$\gcd(\sigma,\sigma+\tau)=1=\gcd(\tau,\sigma+\tau)$ when $\gcd(\sigma,\tau)=1$.
\smallskip

As promised, we show an example where $\Phi_{\sigma,1}\cap\Phi_{\tau,1}
\not= \emptyset$ and $\Phi_{\sigma,2}\cap\Phi_{\tau,1} \not= \emptyset$.
Consider $\phi=70800/1547 < 1000=R$ with $\Delta=\pi/R$. Choose
$\sigma=299$ and $\tau=357$ with $\gcd(\sigma,\tau)=1$.
With 
$$\phi_\sigma=100000/35581, \qquad \ell=68 \le \lceil \sigma/2 \rceil -1 = 149,$$
we have $\phi \in \Phi_{\sigma,1}$.
With 
$$\phi_\tau=6000/1547, \qquad \ell = 81 \le \lceil \tau/2 \rceil-1 = 178,$$ 
we find $\phi \in \Phi_{\tau,1}$.
Unfortunately, since $\phi_\tau = 2R/\sigma-\phi_\sigma$ we also have
$\phi_\tau\in\Phi_{\sigma,2}\cap\Phi_{\tau,1} \not= \emptyset$.
\smallskip

As a last remark, we add that even replacing \eqref{realnyq}
by the stricter constraint
$$|\phi_i|\Delta <\pi/2, \qquad i=1, \ldots, n$$
does not guarantee that each $\phi$ can be identified from only $C_\sigma$ and
$C_\tau$. We illustrate this with a counterexample. Let $\phi=3300/133 <
50=R$ with $\Delta = \pi/(2R)$. With $\sigma=21$ and $\tau=19$ we find
\begin{gather*}
\phi_\sigma=500/133, \qquad (2\pi)/(\sigma\Delta)-\phi_\sigma=2300/399, \\
\phi_\tau=500/133, \qquad (2\pi)/(\tau\Delta)-\phi_\tau=900/133.
\end{gather*}
This leads to $\Phi_{\sigma,1}\cap\Phi_{\tau,1} = \{500/133\}$ and
$\Phi_{\sigma,2}\cap\Phi_{\tau,1} = \{3300/133\}$.


\backmatter



\section*{Declarations}



\subsection*{Funding}
Annie Cuyt and Wen-shin Lee received funding from the European Union's
Horizon 2020 Research and Innovation Staff Exchange program under the MSCA grant agreement No 101008231 (EXPOWER). 
\smallskip

Wen-shin Lee received funding from the Carnegie Trust (Project ``Advancing
exponential analysis: high resolution information from sparse and
regularly sample data''), grant reference RIG009853.

\subsection*{Conflict of interest} 
None.

\subsection*{Ethics approval}
Not applicable

\subsection*{Availability of data and materials}
Code generating the data and running all examples will be available from the
website {\tt cemath.org}.
\subsection*{Authors' contributions}
The publication is the result of joint work, to which both authors contributed equal
efforts.

\end{document}